\documentclass{amsart}
\usepackage{graphicx}
\newtheorem{thm}{Theorem}[section]
\newtheorem{assumption}[thm]{Assumption}
\newtheorem{cor}[thm]{Corollary}
\newtheorem{lem}[thm]{Lemma}
\newtheorem{prop}[thm]{Proposition}
\theoremstyle{definition}

\theoremstyle{remark}
\newtheorem{rem}[thm]{Remark}
\numberwithin{equation}{section}

\newcommand{\eps}{\varepsilon}

\begin{document}

\title{Bandit problems with L\'{e}vy payoff processes}%
\author{Asaf Cohen and Eilon Solan}%
\thanks{The School of Mathematical Sciences, Tel Aviv University, Tel Aviv 69978, Israel. e-mail: cohenasa@post.tau.ac.il}
\thanks{The School of Mathematical Sciences, Tel Aviv University, Tel Aviv 69978, Israel. e-mail: eilons@post.tau.ac.il}

\address{}%
\email{}%

\thanks{We thank Adi Ditkowski and Steven Schochet for their helpful remarks about differential equations, and Sven Rady for useful comments on an earlier version of the paper.}%


\begin{abstract}



We study two-armed L\'{e}vy bandits in continuous-time, which have one safe arm that yields a constant payoff $s$, and one risky arm that can be either of type High or Low; both types yield stochastic payoffs generated by a L\'{e}vy process. The expectation of the L\'{e}vy process when the arm is High is greater than $s$, and lower than s if the arm is Low.

The decision maker (DM) has to choose, at any given time $t$, the fraction of resource to be allocated to each arm over the time interval $[t,t+dt)$.  We show that under proper conditions on the L\'{e}vy processes, there is a unique optimal strategy, which is a cut-off strategy, and we provide an explicit formula for the cut-off and the optimal payoff, as a function of the data of the problem.
We also examine the case where the DM has incorrect prior over the type of the risky arm, and we calculate the expected payoff gained by a DM who plays the optimal strategy that corresponds to the incorrect prior.

In addition, we study two applications of the results: (a) we show how to price information in two-armed L\'{e}vy bandit problem, and (b) we investigate who fares better in two-armed bandit problems: an optimist who assigns to High a probability higher than the true probability, or a pessimist who assigns to High a probability lower than the true probability.
\end{abstract}

\maketitle
\def\high{1}
\def\low{2}

\def\p0{p_0}
\def\p1*{p^{*}}
\def\pp1*{p'}
\def\pt{p_t}
\def\t1eps{t^{\epsilon}_1}
\def\omegapt{( \frac{1-p_t}{p_t})}
\def\omegap{( \frac{1-p}{p})}
\def\fracpt{\frac{p_t}{1-p_t}}
\def\fracpo{\frac{1-p_0}{p_0}}
\def\fracp0{\frac{1-p_0}{p_0}}
\def\fracp00{\frac{p_0}{1-p_0}}
\def\fracp1*{\frac{\p1*}{1-\p1*}}
\def\fraceps{\frac{1-p_0-\epsilon}{p_0+\epsilon}}
\def\fracepseps{\frac{p_0+\epsilon}{1-p_0-\epsilon}}
\def\kt{k_t}
\def\vv{\V^{*}_1}
\def\la{\lambda_1}
\def\laa{\lambda_2}
\def\dYP{dY_{P}^0}
\def\dYPP{dY_{P}^k}
\def\dpit{d \tilde{Y}_B^k}
\def\dpi{d Y_B}
\def\pit{Y_B^1(t)}
\def\pis{Y_B^1(s)}
\def\pita{Y_B^1(t)}

\def\mut{\tilde{\mu}}
\def\Pt{\tilde{P}}
\def\sk{\sqrt{k}}
\def\ta{\theta_{\high}}
\def\tb{\theta_{\low}}
\def\at{\tilde{\mu}_{\high}}
\def\bt{\tilde{\mu}_{\low}}
\def\vb{\bar{\nu}}
\def\vg{\bar{\nu}_{\high}}
\def\vvg{\bar{\nu}_{\low}}
\def\kt{\sqrt{k}}
\def\mat{p \at + (1-p) \bt}
\def\mvg{p \vg + (1-p) \vvg}
\def\mvggag{\hat{p} \vg + (1-\hat{p}) \vvg}

\def\vi{\nu_{\high}}
\def\vvi{\nu_{\low}}
\def\ai{\mu_{\high}}
\def\bi{\mu_{\low}}
\def\gi{g_{\high}}
\def\ggi{g_{\low}}
\def\Hi{H_{\high}}
\def\HHi{H_{\low}}

\def\mvt{p \vg + (1-p) \vvg}
\def\mvh{p \nu_{\high} (dh) + (1-p)\nu_{\low} (dh)}
\def\mvhgag{\hat{p} \nu_{\high} (dh) + (1-\hat{p})\nu_{\low} (dh)}

\def\dZt{ d\tilde{Z}}
\def\vh{\nu_{\high} (dh)}
\def\vvh{\nu_{\low} (dh)}
\def\vhi{\nu_{\high} (dh_j)}
\def\vvhi{\nu_{\low} (dh_j)}
\def\Ph{P_h}

\def\vhvv{\frac{\vh}{\vvh}}
\def\vvhv{\frac{\vvh}{\vh}}
\def\lnvhvv{\ln \left(\frac{ \vhi}{\vvhi}\right)}
\def\lnvvhv{\ln \left(\frac{ \vvhi}{\vhi}\right)}
\def\Fmu{F_{\mu}}

\def\Patau{P_{\theta_{\high}}(\tau<T)}
\def\Patauu{P_{\theta_{\high}}(\tau>T)}
\def\Pbtau{P_{\theta_{\low}}(\tau<T)}
\def\Pbtauu{P_{\theta_{\low}}(\tau>T)}
\def\Pmutau{P_{\theta}(\tau<T)}
\def\Pmutauu{P_{\theta}(\tau>T)}
\def\Pataur{P_{\theta_{\high}}(\tau^r <T)}
\def\Pbtaur{P_{\theta_{\low}}(\tau^r <T)}
\def\Patauur{P_{\theta_{\high}}(\tau^r >T)}
\def\Pbtauur{P_{\theta_{\low}}(\tau^r >T)}

\def\Bta{B^\mu_{\high} (t)}
\def\Btb{B^\mu_{\low} (t)}
\def\Btmu{B^{\mu} (t)}
\def\Bsmu{B^{\mu} (s)}
\def\Bmutau{B^{\mu} (\tau)}

\def\fa{f^{\theta_{\high}}_{(B^{\mu_{\high}}(\tau),\tau)| \tau < T} (x,t)}
\def\fb{f^{\theta_{\low}}_{(B^{\mu_{\low}}(\tau),\tau)| \tau < T} (x,t)}

\def\fmu{f^{\theta}_{(B^{\mu}(\tau),\tau)| \tau < T} (x,t)}

\def\fmu{f^{\theta_i}_{(B^{\mu_i}(\tau),\tau)| \tau < T} (x,t)}

\def\Xkn{N^n_{\kappa_n}}
\def\Ykn{Y^n_{\kappa_n}}
\def\Zkn{Z_{\kappa_n}}

\def\Fns{\mathcal{F}^n_s}
\def\Gns{\mathcal{G}^n_s}
\def\Fnt{\mathcal{F}^n_t}
\def\Gnt{\mathcal{G}^n_t}
\def\Fms{\mathcal{F}^m_s}
\def\Gms{\mathcal{G}^m_s}
\def\Fmt{\mathcal{F}^m_t}
\def\Gms{\mathcal{G}^m_t}

\def\Fn{\mathcal{F}^n}
\def\Fm{\mathcal{F}^m}
\def\Gn{\mathcal{G}^n}
\def\Gm{\mathcal{G}^m}

\def\F{\mathcal{F}}
\def\G{\mathcal{G}}

\def\Ft{\mathcal{F}_t}
\def\Gt{\mathcal{G}_t}

\def\Fs{\mathcal{F}_s}
\def\Gs{\mathcal{G}_s}

\def\Fl{\mathcal{F}^l}
\def\Fm{\mathcal{F}^m}

\def\eps{\epsilon}

\def\Lknn{L^n_{\kappa_n}}
\def\Nknn{N^n_{\kappa_n}}

\def\Mk{M_{\kappa}}
\def\Lkn{L^n_{\kappa}}
\def\Nkn{N^n_{\kappa}}

\def\Nknn{N^n_{\kappa_n}}

\def\Mkn{M_{\kappa_n}}

\def\ata{\tilde{\mu}^a_{\high}}
\def\bta{\tilde{\mu}^a_{\low}}
\def\vba{\bar{\nu}^a}
\def\vga{\bar{\nu}^a_{\high}}
\def\vvga{\bar{\nu}^a_{\low}}

\def\atb{\tilde{\mu}^b_{\high}}
\def\btb{\tilde{\mu}^b_{\low}}
\def\vbb{\bar{\nu}^b}
\def\vgb{\bar{\nu}^b_{\high}}
\def\vvgb{\bar{\nu}^b_{\low}}

\def\atc{\tilde{\mu}^c_{\high}}
\def\btc{\tilde{\mu}^c_{\low}}
\def\vbc{\bar{\nu}^c}
\def\vgc{\bar{\nu}^c_{\high}}
\def\vvgc{\bar{\nu}^c_{\low}}

\def\via{\nu^a_{\high}}
\def\vvia{\nu^a_{\low}}
\def\aia{\mu^a_{\high}}
\def\bia{\mu^a_{\low}}
\def\gia{g^a_{\high}}
\def\ggia{g^a_{\low}}
\def\Hia{H^a_{\high}}
\def\HHia{H^a_{\low}}

\def\vib{\nu^b_{\high}}
\def\vvib{\nu^b_{\low}}
\def\aib{\mu^b_{\high}}
\def\bib{\mu^b_{\low}}
\def\gib{g^b_{\high}}
\def\ggib{g^b_{\low}}
\def\Hib{H^b_{\high}}
\def\HHib{H^b_{\low}}

\def\vic{\nu^c_{\high}}
\def\vvic{\nu^c_{\low}}
\def\aic{\mu^c_{\high}}
\def\bic{\mu^c_{\low}}
\def\gic{g^c_{\high}}
\def\ggic{g^c_{\low}}
\def\Hic{H^c_{\high}}
\def\HHic{H^c_{\low}}

\def\vha{\nu^a_{\high} (dh^a)}
\def\vvha{\nu^a_{\low} (dh^a)}
\def\vhia{\nu^a_{\high} (dh^a_j)}
\def\vvhia{\nu^a_{\low} (dh^a_j)}

\def\vhvva{\frac{\vha}{\vvha}}
\def\vvhva{\frac{\vvha}{\vha}}
\def\lnvhvva{\ln \left(\frac{ \vhia}{\vvhia}\right)}
\def\lnvvhva{\ln \left(\frac{ \vvhia}{\vhia}\right)}

\def\vhb{\nu^b_{\high} (dh^b)}
\def\vvhb{\nu^b_{\low} (dh^b)}
\def\vhib{\nu^b_{\high} (dh^b_j)}
\def\vvhib{\nu^b_{\low} (dh^b_j)}

\def\vhvvb{\frac{\vhb}{\vvhb}}
\def\vvhvb{\frac{\vvhb}{\vhb}}
\def\lnvhvvb{\ln \left(\frac{ \vhib}{\vvhib}\right)}
\def\lnvvhvb{\ln \left(\frac{ \vvhib}{\vhib}\right)}

\def\vhc{\nu^c_{\high} (dh^c)}
\def\vvhc{\nu^c_{\low} (dh^c)}
\def\vhic{\nu^c_{\high} (dh^c_j)}
\def\vvhic{\nu^c_{\low} (dh^c_j)}

\def\vhvvc{\frac{\vhc}{\vvhc}}
\def\vvhvc{\frac{\vvhc}{\vhc}}
\def\lnvhvvc{\ln \left(\frac{ \vhic}{\vvhic}\right)}
\def\lnvvhvc{\ln \left(\frac{ \vvhic}{\vhic}\right)}

\def\vgj{\bar{\nu}^j_{\high}}
\def\vvgj{\bar{\nu}^j_{\low}}

\def\muga{\tilde{\mu}^a}
\def\mugb{\tilde{\mu}^b}

\def\faj{f^{\theta_{\high}}_{(B^{\muga,\mugb}(\tau),\tau)| \tau < T} (x,t)}
\def\fbj{f^{\theta_{\low}}_{(B^{\muga,\mugb}(\tau),\tau)| \tau < T} (x,t)}

\section{Introduction}

Consider a firm that has to determine,
on an ongoing basis,
how much to invest in the research of new technologies for its next line of products.
The firm faces a tradeoff between exploration and exploitation:
on the one hand, it can adopt the technology that seems most successful according to the research conducted so far,
thereby exploiting its investment in research,
but on the other hand, it could continue investing in various technologies,
in the hope of finding an even better technology for its products.
If the firm decides to stop investing in a given technology,
then no information will be obtained on that technology,
so even if it is actually better than the finally adopted technology, it will never be adopted.

A similar tradeoff between exploration and exploitation arises, e.g., in the market of venture capital funds,
where each fund has to decide in which start-up companies to invest,
and in clinical trials, where pharmaceutical companies have to decide which new drugs or treatments to explore.

To concentrate on the trade off between exploration and exploitation,
one assumes that there are no exogenous factors that affect the
firm's decision (such as new technologies or drugs that are
introduced by competitors). The optimization problem that the firm
faces has been modeled in the literature as a multi-arm bandit
problem (see, e.g., Rothschild (1974), Bergemann and V\"{a}lim\"{a}ki (2006), Keller, Rady
and Cripps (2005), Besanko and Wu (2008), Klein and Rady (2008), Moscarini and Squintani (2004), Roberts and Weitzman (1981), Weitzman (1979)): a decision maker (DM) has
finitely many actions, called arms, each one yields a payoff with
an unknown distribution, that is taken from a finite set of
distributions. Each time the DM chooses an arm, he
obtains a payoff, and improves his information regarding the
correct payoff distribution of the arm he has just chosen.

Gittins and Jones (1979)  
proved that, in discrete-time, the optimal strategy of the DM has a particularly simple form: at every period the DM
calculates a real number, an index, for each arm, based on
past observations of that arm, and he chooses
the arm with the highest index. It turns out that
to calculate the index of an arm it is sufficient to consider an
auxiliary problem with two arms: the arm for which we calculate the
index, and an arm that yields a constant payoff. The literature
therefore focuses on such problems, called two-armed bandit
problems.

Once the optimality of the index strategy is guaranteed,
one looks for the relation between that data of the game and the index.
Explicit formulas for the index
when the payoff is one of two distributions that have a simple form has been established in the literature.
Berry and Friestedt (1985) provide the solution to the problem in discrete-time,
e.g., when the payoff distribution is one of two Bernoulli distributions,
and in continuous-time,
e.g., when the payoff distribution is one of two Brownian motions.
In continuous-time, by studying the dynamic programming
equation that describes the problem,
Keller, Rady and Cripps (2005)
and Keller and Rady (2008)
provided an explicit form for the index when the payoff's distribution
is Poisson.%
\footnote{These authors also studied the strategic setup, in which
several DMs have the same set of arms, and their arms'
payoff distributions are the same (and unknown), and they compared
the cooperative solution to the non-cooperative solution.}
When the payoff distribution is known, Karatzas (1984) characterized the index when the payoff's distribution is a diffusion process,
and Kaspi and Mandelbaum (1995)
characterized the index when the payoff's distribution is a L\'{e}vy process,
and they obtained an explicit form for the index for special distributions.

In practice, payoff processes have a complex form,
exhibiting both small random changes, that can be modeled by a Brownian motion,
and large shocks that can be modeled as arriving at a Poisson rate. A stochastic process that incorporates these two types of changes is the L\'{e}vy process.
In fact,
Carr and Wu (2004) argue that almost all economic phenomena can be described by time shifts of L\'{e}vy processes.
Therefore, it is desirable to study the bandit problem
where the payoff distribution is one of finitely many L\'{e}vy distributions.

In the present paper we provide an explicit solution to the two-armed bandit problem
where the payoff distribution is one of two L\'{e}vy processes.
We assume
that one distribution, called {\em High}, dominates the other,
called {\em Low}, in a strong sense (see Assumption \ref{assumption} below). To
eliminate trivial cases, we assume that the expected payoff that
is generated by the safe arm is lower than the expected payoff
generated by the High distribution, and higher than the expected
payoff generated by the Low distribution.

In such a case in discrete-time, the optimal strategy is a cut-off
strategy: the DM keeps on experimenting as long as the
posterior belief that the distribution is High is higher than some
cut-off point, and, once the posterior probability that the distribution
is High falls below the cut-off point, the DM switches to
the safe arm. 
We prove that when the two payoff
distributions are L\'{e}vy processes that satisfy several requirements,
the optimal strategy is a cut-off strategy, and we provide an explicit
expression for the cut-off point, in terms of the data of the problem.
When particularized to the models studied by Kaspi and Mandelbaum
(1995), Bolton and Harris (1999), Keller, Rady and Cripps (2005) and Keller and Rady (2008),
our expression reduces to the expressions that they obtained.

Apart of unifying previous results, our characterization shows that the special form of the optimal payoff derived by Bolton and Harris (1999) and Keller, Rady and Cripps (2005) is valid in a general setup: the optimal payoff is the sum of the expected payoff if no information is available, and of an option value, that measures the expected gain from the ability to experiment. It also shows that the data of the problem can be divided into information-relevant parameters and payoff-relevant parameters; the information-relevant parameters can be summarized in a single real number, and the payoff-relevant parameters are the expectations of the processes that contribute to the DM's payoff.
Finally, the characterization allows one to derive comparative statics on the optimal cut-off and payoff. For example, as the discount rate increases, or the signals become less informative, the cut-off point increases but the DM's optimal payoff decreases.

It is often the case that the DM holds an incorrect prior,
and plays optimally given that prior (for further discussion and literature review, see Section \ref{incorrect prior}).
We provide an explicit expression for the value function in this case.
Our technique provides a new description of the optimal strategy with a time dependent cut-off.

So far we have assumed that all the information that the DM has is the payoff process.
Sometimes, the DM has information that does not contribute to the payoff process,
yet it helps him learn the type of risky arm.
For example, scientific discoveries made by other firms in other markets may shed light
on the appropriateness of a given technology to a product that a firm develops.
In addition, part of the payoff of the DM may not be observed by the DM.

In Section \ref{information} we study a bandit problem, in which the risky arm generates three L\'{e}vy processes - the first is observed by the DM and contributes to the payoff, the second is observed by the DM and does not contribute to the payoff, and the third is not observed by the DM but it contributes to the payoff. We provide an explicit expression for the cut-off and for the optimal payoff of the DM. That generalizes the expression we found when the risky arm generates only one process. This analysis clarifies the distinction between information-relevant data and payoff-relevant data.

%
We conclude the paper by applying our characterization to compare the effects of optimism and pessimism in bandit problems.
A DM is called optimist if his prior probability that the payoff's distribution is High
is higher than the true probability,
and he is called a pessimist if his prior probability that the payoff's distribution is High
is lower than the true probability.
Using our characterization we find 
that unless the pessimist assigns high probability to the High type and the two DM's are sufficiently patient, an optimist will fare better than a pessimist.


The rest of the paper is arranged as follows.
The model and the main results appear in Section 2. Directions for future research appear in Section 3. All proofs appear in Section 4.

\newpage
\section{the model and the main results}

\subsection{Reminder about L\'{e}vy Processes}\label{s:reminder}
 L\'{e}vy processes are the continuous-time analog of discrete-time random walks with i.i.d. increments. A \emph{L\'{e}vy process} $X=(X(t))_{t\geq 0}$ is a continuous-time stochastic process that (a) starts at the origin: $X(0)=0$, (b) admits c\`{a}dl\`{a}g modification,\footnote{That is, it is continuous from the right, and has limits from the left: for every $t_0$, the limit $X(t_0 -):=\underset{t\nearrow t_0}{\lim}X(t)$ exists a.s. and $X(t_0) =\underset{t\searrow t_0}{\lim}X(t)$. } and (c) has stationary independent increments.
See figures 1 and 2 for a generic path of a L\'{e}vy process. A few examples of L\'{e}vy processes are a Brownian motion, a Poisson process, and a compound Poisson process. The latter is a continuous-time process in which jumps arrive according to a Poisson process and the jumps are i.i.d.\footnote{Formally, let $\lambda>0$, and let $D$ be a distribution over $\mathbb{R}\backslash\{0\}$. A \emph{compound Poisson process} with rate $\lambda$ and jump size distribution $D$ is a continuous-time stochastic process given by $X(t) = \overset{N(t)}{\underset{i=1} {\sum}} D_i$, where $N(t) $ is a Poisson process with rate $\lambda$, and $D_i$ are i.i.d. random variables, with distribution function $D$, which are also independent of $(N(t))_{t\geq 0}$.}




We now present the L\'{e}vy-Ito decomposition of L\'{e}vy processes. Let $(X(t))$ be a L\'{e}vy process. For every Borel measurable set $A\subseteq \mathbb{R}\backslash \{0\}$, define:
\[
\nu (A):= E\left[\sharp\{0\leq s \leq 1 | \Delta X(t): = X(t) -X(t-) \in A\}\right].
\]
This is the expected number of jumps with size in $A$ that occurs up to time $1$. One can verify that $\nu$ is a measure, called the \emph{L\'{e}vy measure} of the process. The quantity $\nu(\mathbb{R}\backslash\{0\})$ is the expected number of jumps that occur up to time $1$. If $\nu(\mathbb{R}\backslash\{0\})=\infty$, then the number of jumps in the time interval $[0,1]$, and therefore in any compact time interval, is infinite a.s., and we say that the L\'{e}vy measure is infinite. If $\nu (\mathbb{R}\backslash\{0\}) < \infty,$ then the expected number of jumps in any compact time interval is finite a.s., and we say that the L\'{e}vy measure is finite. In this paper we study L\'{e}vy processes that have finite L\'{e}vy measures. The results can be generalized for L\'{e}vy processes with infinite L\'{e}vy measures, see Cohen and Solan (2009).

The L\'{e}vy-Ito decomposition (see Applebaum (2004)) states that every L\'{e}vy process with finite L\'{e}vy measure can be represented as follows:
\begin{equation}\label{levy-ito}
X(t) = \mu t +\sigma Z(t) + L_{\nu}(t)   ,
\end{equation}
where $\mu t$ is a linear drift, $\sigma Z(t)$ is a Brownian motion with standard deviation $\sigma$, and  $L_{\nu}(t)$ is a compound Poisson process with L\'{e}vy measure $\nu$ which is independent of $\sigma Z(t)$: jumps arrive at a Poisson rate with expectation $\nu (\mathbb{R}\backslash\{0\})$, and the distribution of each jump is given by the distribution function $\frac{\nu(dh)}{\nu(\mathbb{R}\backslash\{0\})}$.

\subsection{L\'{e}vy Bandits}\label{s:finite}
A DM operates a two-armed bandit machine in continuous-time, with a safe arm that yields a constant payoff $s$, and a risky arm that yields a stochastic payoff $(X(t))$. The risky arm can be of two types, High or Low. We denote the arm's type by $\theta$: if the type is High (resp. Low) we set $\theta = \ta$ (resp. $\tb$). If $\theta =\theta_i$, $i\in\{\high,\low\}$, the risky arm yields payoff $(X_i (t)),$ which is a L\'{e}vy process.
We assume throughout that the L\'{e}vy measures of both $ (X_{\high} (t))$ and $(X_{\low} (t))$ are finite, and therefore a.s. there are only finitely many jumps in each compact time interval. Denote the L\'{e}vy-Ito decomposition of $(X_i (t))$, $i\in\{\high,\low\},$ by   $X_i(t) = \mu_i t + \sigma_i Z(t) + L_{\nu_i}(t)$.



Set $\bar{\nu_i} := \nu_i (\mathbb{R}\setminus\{0\})$, and denote by $H_i:= \int h\nu_i (dh)/\bar{\nu_i }$ the expected jump size of $(X_i (t))$.\footnote{In order to simplify notation, we denote $\int h\nu_i (dh):=\underset{\mathbb{R}\backslash\{0\}}{\int} h\nu_i (dh)$.} We assume that $H_i$ is finite. The quantity $\int h\nu_i (dh)= \bar{\nu_i } H_i$ is the contribution of the compound Poisson process to the instantaneous payoff. The expectation of the risky arm at time $t=1$ if $\theta=\theta_i$ is $g_i :=E[X_i(1)] =  \bar{\nu_i} H_i +\mu_i$.

Throughout we make the following assumption, which states that the High type is better than the Low type in a strong sense.

\begin{assumption}\label{assumption}$\\$
A1. $-\infty<\ggi<s<\gi<\infty$.\\
A2. $\sigma_{\high}=\sigma_{\low}$.\\
A3. For every $A\in\mathcal{B}(\mathbb{R}\setminus \{0\}),\; \; \vvi (A)\leq\vi (A)<\infty$.
\end{assumption}


Assumptions A1 and A2 rule out trivial cases.
Assumption A1 merely says that the High (resp. Low) type provides higher (resp. lower) expected payoff than the safe arm. Assumption A2 states that the Brownian motion component of both the High type and Low type have the same standard deviation. Otherwise, since the realized path reveals the standard deviation, the DM can distinguish between the arms in any infinitesimal time interval.

The third part of the assumption is less innocuous; it requires 
that the L\'{e}vy measure of the High type will dominate the L\'{e}vy measure of the Low type in a strong sense:
roughly, jumps of any size $h$ occur more often (or at the same rate) under the high type than under the Low type.

A consequence of this assumption is that jumps always provide good news, and (weakly) increase the posterior probability of the High type.

At each time instance $t$, the DM chooses the proportion of time to devote to each of the two arms. If he chooses to devote a proportion $k$ of the current time instance to the risky arm, then the instantaneous payoff $dY^k$ 
%
is the sum of several terms:
\begin{itemize}
\item
$(1-k)sdt$, which is the contribution of the safe arm;
\item
$k\mu_i dt$, which is the contribution of the linear drift;
\item
$\sqrt{k}\sigma dZ(t)$, which is the contribution of the Brownian motion;\footnote{By devoting the proportion $\kappa$ of the time interval $[t,t+dt)$ to the risky arm, the variance of the continuous part of the payoff is $\kappa\sigma^2 dt$. This explains the scaling parameter $\sqrt{k}\sigma$. For a qualitative explanation for this form, see Bolton and Harris (1999).}
\item
$kH_i dt$, which is the contribution of the compound Poisson process.
\end{itemize}


A strategy $\kappa$ is a (measurable) function, that assigns to each history a number in the interval $[0,1]$, that is interpreted as the amount of time in the interval $[t,t+dt)$ devoted to the risky arm.
In continuous time, it is usually assumed that a strategy is predictable, that is, to determine the behavior at time $t$, it is sufficient to know the history strictly before time $t$. Formally, $\kappa_t$ is $\F_{t-}^I$-measurable, where
$I(t) := \int_0^t re^{-rt} dY^{\kappa} (t)$ is the stochastic discounted payoff using the strategy $\kappa$, and
$\F_{t-}^I$ is the $\sigma$-algebra generated by the stochastic process $(I(t))$ of the discounted payoff with discount rate $r$, up to (excluding) time $t$.

As is well known , in continuous-time the play path need not be uniquely defined, and therefore one should be careful in defining the set of strategies available to the DM, and the play path that a strategy defines. We circumvent this issue by arguing that an optimal strategy must solve a certain Functional Differential Equation (FDE), by showing that this FDE has a unique solution, and by exhibiting this solution. As for discrete-time, this solution turns out to be a cut-off strategy, which we now define.

Let $p_t:=P(\theta=\ta |\F_{t-}^I)$ be the posterior belief at time $t$ that the risky arm is High. A strategy $\kappa$ is a \emph{cut-off strategy} with cut-off point $p$ if the DM chooses the safe arm (with probability 1) whenever $p_t\leq p $, and the risky arm (with probability 1) whenever $p_t>p$. We now argue that Assumption A3 guarantees that the play under a cut-off strategy is well defined.
Suppose that $p_0>p $. Then the DM chooses the risky arm until the first time $t$ that satisfies $p_t=p$. Since the L\'{e}vy payoff processes have finite measures, Assumption A3 implies that in an infinitesimal interval after time $t$ the posterior belief will drop below the cut-off. Indeed,
the first time $t$ satisfying $p_t=p$ is a predictable stopping time, and therefore, the probability that a jump occurs atthat time is zero (see Bertoin (1996)). Therefore, $P(\theta=\ta |\F_{t}^I)=P(\theta=\ta |\F_{t-}^I)=p$.
If there is a Brownian motion component,  its fluctuations will cause the posterior to drop below $p$ in an infinitesimal time interval after time $t$. If there is no Brownian motion component, since the L\'{e}vy measure is finite, there are no jumps in an infinitesimal time interval, and by Assumption A3 the compound Poisson process will cause the posterior to drop below $p$. This implies that the play under a cut-off strategy is well defined.

\subsection{The Optimal Strategy}\label{s:optimal}

The expected discounted payoff under a strategy $\kappa$ when the prior is $p_0 =p$ is
\begin{align}\notag
V_{\kappa}(p) &= E \left[\int_0^{\infty} r e^{-rt}dY^{\kappa} (t)\right] \\\notag
              &= p E \left[\left.\int_0^{\infty} r e^{-rt}dY^{\kappa} (t)\right|\ta\right]+(1-p) E \left[\left.\int_0^{\infty} r e^{-rt}dY^{\kappa} (t)\right|\tb\right].\\\notag
\end{align}
Let $U(p) = \underset{\kappa}{sup}V_{\kappa}(p)$ be the maximal payoff the DM can achieve. As we show below, the DM has an optimal strategy, so in fact the supremum in the definition of $U(p)$ is achieved. The function $ V_{\kappa}(p)$ is linear with respect to $p$, and therefore $U(p)$, as the supremum of linear functions, is convex. By always choosing the safe arm, the DM can achieve at least $s$; since $U(0)=s$, the convexity of $U(p)$ implies that $U$ is non-decreasing.

\begin{prop}\label{convex}
$U(p)$ is monotone non-decreasing, convex, and continuous in $p$.
\end{prop}
It follows from Proposition \ref{convex} that there is $p^*$ such that $U(p)=s$ if $p\leq p^*$ and $U(p) >s$ otherwise, so that the strategy $\kappa\equiv 0$ that always chooses the safe arm is optimal for prior beliefs in $[0,\p1*]$.\\\\
Our first theorem states that 
there is a unique optimal strategy, which is a cut-off strategy. Moreover it provides the exact cut-off point and the corresponding expected payoff in terms of the data of the problem. Let $\alpha$ be the unique solution of
\begin{equation}\label{eta}
f(\eta) :=  \int\vvh\left(\frac{\vvh}{\vh}\right)^{\eta} + \eta (\vg-\vvg)-\vvg +\frac{1}{2}(\eta+1)\eta \left(\dfrac{\ai-\bi}{\sigma}\right)^2 - r = 0
\end{equation}
in $(0,\infty)$.
The existence and uniqueness of such a solution are proved in Lemma \ref{alpha} below.
Observe that $\frac{\vvh}{\vh}$, the Radon-Nikodym derivative, exists by Assumption A3.\footnote{To ensure the existence of the Radon-Nikodym derivative, and therefore the form of the solution, one does not need the full power of Assumption A3. Its full power will be used within the proof.} 

\begin{thm}\label{Up1}
Denote $\p1*:= \frac{\alpha (s-\ggi)}{(\alpha+1)(\gi-s) + \alpha(s-\ggi)}$. Under Assumptions A1-A3, the unique optimal strategy is\\
$\kappa^* =
\begin{cases}
0               &\text{if $p \leq \p1*$}, \\
1               &\text{if $p > \p1*$}.
\end{cases}
$\\
The expected payoff under $\kappa^*$ is
\begin{equation}\label{E:Up1}
U(p) = V_{\kappa^*} (p) =
\begin{cases}
s               &\text{if $p \leq \p1*$}, \\
p\gi +(1-p)\ggi +C_\alpha (1-p)(\frac{1-p}{p})^{\alpha}               &\text{if $p > \p1*$},
\end{cases}
\end{equation}
$\\$
where $C_\alpha = \frac{s-\ggi - \p1* (\gi -\ggi)}{(1-\p1*) \left( \frac{1-\p1*}{\p1*}\right)^{\alpha}}.$\\
\end{thm}

The term $p\gi +(1-p)\ggi$ in ~(\ref{E:Up1}) is the expected payoff for a DM choosing only the risky arm. Thus, $C_\alpha (1-p)(\frac{1-p}{p})^{\alpha}$ is the option value for the ability to switch to the safe arm.
The quantity $\alpha$ summarizes all the information-relevant parameters of the problem (see also Section \ref{information}). Apart from that, the parameters of the payoff processes that determine the cut-off point $\p1*$ and the optimal payoff $U(p)$ are the expected payoffs $\gi$ and $\ggi$.

The function in ~(\ref{E:Up1}) has the same structure as the solution of Bolton and Harris (1999) and Keller, Rady and Cripps (2005) for the one agent problem.
In  Bolton and Harris (1999), the only component in the risky arm is the Brownian motion with drift. Therefore $\nu_i \equiv 0$, so that $\gi  = \ai,\; \ggi = \bi,$ and $ \alpha = (-1+\sqrt{1+8r\sigma^2 / (\ai -\bi)^2})/2$.
In Keller, Rady and Cripps (2005), the risky arm is either the constant zero (Low type, so that $\vvi \equiv 0$), or yields a payoff $\bar{h}$ according to a Poisson process with rate $\lambda$ (High type). If the risky arm is High, the only component in the L\'{e}vy-Ito decomposition is the compound Poisson component, and $\vi (\bar{h})  = \lambda$ and zero otherwise.
Therefore,  $\gi = \lambda \bar{h}, \; \ggi = 0,$ and $\alpha  = r/\lambda$.\\


The explicit form of the cut-off point $\p1*$ and of the value function $U$ allows us to derive simple comparative statics.
As is well known, a DM who plays optimally switches to the safe arm later than a myopic DM, and indeed, $\p1*$ is smaller than the myopic cut-off point $p^m:=\frac{s-\ggi}{\gi-\ggi}$.
Furthermore, the cut-off point $\p1*$ is an increasing function of $\alpha$. As can be expected, $\alpha$ (and therefore also $\p1*$) increases at the discount rate $r$ and at $\vvi (dh)$, and it decreases at $\vi (dh)$ and at $|\ai-\bi|$: the DM switches to the safe arm earlier as the discount rate increases, as jumps provide less information, or as the difference between the drifts of the two types increases.\footnote{Moreover, $\alpha(r=0) = 0$ and $\alpha(r=\infty) = \infty$. }
Furthermore, as long as $p>\p1*$, the value function $U(p)$ decreases in $\alpha$. Thus, decreasing the discount rate, increasing the informativeness of the jumps and the difference between the drifts is beneficial to the DM.

In Section \ref{incorrect prior} we extend the results to the case that the prior belief of the DM is not the true prior $p_0$.
In Sections \ref{s:optimism} and \ref{information} we provide two applications to our techniques and results.

\begin{rem} \textbf{Incomplete information of $p_0$.}
Suppose that the DM does not know the prior belief $p_0$, but rather has some belief $\varphi$ over $p_0$. That is, $p_0$ is chosen at the outset according to $\varphi$, and is not told to the DM. From the DM's point of view, the situation is equivalent to an auxiliary problem in which the prior probability of the High type is the expectation of the corresponding probability in the original problem, and therefore Theorem \ref{Up1} provides the optimal strategy in this case as well.
\end{rem}



\subsection{The Payoff with Incorrect Prior}\label{incorrect prior}

In decision problems it is usually assumed that the DM
either knows the true state of nature,
or has some prior distribution over the set of states of nature.
Experiments show that the prior distribution that DMs have
is often different than the true prior.
The phenomenon of overconfidence --
assigning too high a probability to the good state of nature --
has been observed in various areas
(Svenson (1981), Baumhart (1968), Larwood and Whittaker (1977), Cross (1977),
Weinstein (1980), Camerer and Lovallo (1999)).
Babcock and Loewenstein (1997) argue that biases in bargaining may be self-serving,
and Heifetz, Shannon and Spiegel (2007) show that biases of preferences may be stable.

In every decision problem,
a DM who correctly perceives the prior distribution will fare better than a DM who has some bias,
and believes that the prior distribution is different than the correct one.
Indeed, an optimal strategy of a DM who correctly perceives the prior distribution
is a strategy that yields the highest possible gains for this prior distribution,
so it yields at least as much as any other strategy,
in particular, optimal strategies for incorrect prior distributions.


Denote the initial belief of the DM for the High type by $q_0$, and suppose that it may be different from the true probability $p_0$. 
By Theorem \ref{Up1}, the optimal strategy of the DM is a cut-off strategy, however, since he has an incorrect prior, he does not switch to the safe arm at the optimal time.
In this section we give an exact formula for the payoff, assuming the DM plays optimally given his belief.
We will also describe the optimal strategy from a different point of view, not as a cut-off strategy. This point of view is arguably closer to the way people perceive the decision problem that the DM faces.

Suppose that until time $t$, the DM chose the risky arm, and observed the jumps $h_1,...h_n$ from the compound Poisson component of the payoff process. Let $Y_B^1 (t)$ be the Brownian motion with drift component of the payoff process from the risky arm at time $t$. Note that $Y_B^1 (t) \sim N(\mu t  , \sigma^2 t )$. The posterior belief $q_t: = P_t(\ta |h_1,...h_n; \pit; q_0)$ of the DM is:\footnote{$\prod_{t-} \nu_i (dh_j)$ is the product of $\nu_i (dh_j)$ over all jumps $h_j$ that occur up to time $t-$. Similarly, we use the notation $\sum_{t-}$ for the sum over all jumps up to time $t-$. }
\begin{align}\label{posterior}
  q_t &= \frac{q_0 \frac{1}{\sqrt{2\pi t}\sigma} e^{-\frac{(\pit-\ai t)^2}{2\sigma^2 t}} e^{-\vg t}\prod_{t-}\vhi }{q_0 \frac{1}{\sqrt{2\pi t}\sigma} e^{-\frac{(\pit-\ai t)^2}{2\sigma^2 t}} e^{-\vg t} \prod_{t-}\vhi+  (1- q_0) \frac{1}{\sqrt{2\pi t}\sigma} e^{-\frac{(\pit-\bi t)^2}{2\sigma^2 t}} e^{-\vvg t}\prod_{t-}\vvhi} \\\notag\\\notag & = \frac{q_0 e^{\ai \pit /\sigma^2 - \ai^2 t /2\sigma^2 } e^{-\vg t} \prod_{t-} \vhi}{ q_0 e^{\ai \pit /\sigma^2 - \ai^2 t /2\sigma^2 } e^{-\vg t} \prod_{t-} \vhi  +   (1- q_0) e^{\bi \pit /\sigma^2 - \bi^2 t /2\sigma^2 } e^{-\vvg t} \prod_{t-} \vvhi}.\\\notag
\end{align}
Indeed, $\frac{1}{\sqrt{2 \pi t}\sigma} e^{-\frac{(\pit-\mu_i t)^2}{2\sigma^2 t}}$ is the probability of receiving the payoff $Y_B^1 (t)$, given the type $\theta_i$, and $e^{-\bar{\nu}_i t}\frac{(\bar{\nu}_i t)^n}{n!}\prod_{t-}\frac{\nu_i (dh_j)}{\bar{\nu}_i}  $ is the probability of receiving the $n$ jumps that occurred until time $t$, given the type $\theta_i$. The first equality in ~(\ref{posterior}) is the Bayesian belief updating, using the independence of the components in the L\'{e}vy-Ito decomposition, given the type of the risky arm, and the second equality is obtained by eliminating common components. For a generic path of the process $(q_t)$, see figure 3.\\\\

Suppose the DM 
follows a cut-off strategy  $\kappa'$ with cut-off point $p'$:\\
$\kappa' =
\begin{cases}
0               &\text{if $p \leq p'$}, \\
1               &\text{if $p > p'$}.
\end{cases}
$\\

If $q_0 \leq p'$, the DM will always choose the safe arm. If $q_0>p'$, the DM will initially choose the risky arm. 
The DM chooses the risky arm as long as $q_t > p'$, which, by Eq.\ ~(\ref{posterior}), is equivalent to:
\begin{equation}
\dfrac{q_0 e^{\ai \pit /\sigma^2 - \ai^2 t /2\sigma^2 -\vg t}    }{(1-q_0) e^{\bi \pit /\sigma^2 - \bi^2 t /2\sigma^2 -\vvg t}  } \prod_{t-} \vhvv > \dfrac{p'}{1-p'}.
\end{equation}
Without loss of generality assume that $\ai-\bi > 0$; the case $\ai-\bi<0$ is handled similarly, and provides the same results. By taking the natural logarithm, and rearranging the resulting terms, we obtain that this inequality is equivalent to:
%
\begin{align}\label{posteriorq}
\frac{1}{\sigma}\pit &> \left[\left( \frac{\ai + \bi}{2\sigma}\right) + \frac{\sigma(\vg-\vvg)}{\ai-\bi}\right]t  - \frac{\sigma}{\ai -\bi}\times \left[\ln\left(\frac{q_0}{1-q_0}\right)-\ln\left(\frac{p'}{1-p'}\right)\right]\\\notag
&- \frac{\sigma}{\ai -\bi}\sum_{t-} \lnvhvv.
\end{align}
%
%
The right-hand side in ~(\ref{posteriorq}) is a piecewise linear function $F\cdot t-E-G_t$ of $t$, where the slope $F:=  \left( \frac{\ai + \bi}{2\sigma}\right) + \frac{\sigma(\vg-\vvg)}{\ai-\bi}$ is independent of $t$, the intercept at $t=0$ is $E : = \frac{\sigma}{\ai -\bi}\times \left[\ln\left(\frac{q_0}{1-q_0}\right)-\ln\left(\frac{p'}{1-p'}\right)\right]$,
and $G_t: = \frac{\sigma}{\ai -\bi}\sum_{t-} \lnvhvv $.
Denote by $G_h :=\frac{\sigma}{\ai -\bi}\ln\left(\vhvv\right) $ the contribution of a jump of size $h$ to the intercept.\\\\

From Eq.\ ~(\ref{posteriorq}) we obtain the following alternative description of the optimal strategy:
The DM has a time-dependent cut-off which is piecewise linear. The slope of the cut-off function is always $F$, and whenever there is a jump of size $h$, the cut-off decreases by $G_h$ (see Figure 4). The DM chooses the risky arm as long as his current payoff from the continuous part of the L\'{e}vy process, $\pit$, exceeds the cut-off.

Thus, at first the DM plays until the payoff from the continuous part divided by the standard deviation satisfies: $\frac{1}{\sigma} \pit  \leq F\cdot t  -E$; if a jump of size $h$ occurs before he switches to the safe arm, then the intercept decreases by $G_h$, and this behavior repeats itself.
If there is no Brownian motion component, then $\pit \equiv 0$ and $F=0$; the DM chooses the risky arm for a fixed amount of time, and then switches to the safe arm, unless a jump occurs before the switch; if a jump occurred, the amount of time to choose the risky arm increases, as a function of $ \frac{\vi (dh)}{\vvi (dh)}$.


\begin{thm}\label{UpE1}
Under Assumption \ref{assumption}, for every $p_0\in [0,1]$, the payoff $U(p_0,q_0)$ of a DM who uses a cut-off strategy with cut-off point $p'$, is as follows:\footnote{We omit the dependence of $U$ on $p'$, $s$, $g_1$, $g_2$. Recall that $\alpha$ is the unique solution of the equation $f(\eta)=0$ (see Eq.\ ~(\ref{alpha})).} if $q_0\geq p'$,
\begin{align}\label{E:UpE1}
U(p_0,q_0) &=
(s-\gi)p_0 \left(  \frac{1-q_0}{q_0}   \right)^{\alpha + 1} \left(  \frac{p'}{1-p'}   \right)^{\alpha + 1}
+   (s-\ggi)(1-p_0) \left(  \frac{1-q_0}{q_0}   \right)^{\alpha } \left(  \frac{p'}{1-p'}   \right)^{\alpha }\\\notag
&+\gi p_0 + \ggi (1-p_0),
\end{align}
while $U(p_0,q_0) = s$ if $q_0<p'$.
\end{thm}
One can verify that when the DM holds the correct prior, and plays according to the optimal strategy $\kappa^*$, Eq.\ ~(\ref{E:UpE1}) coincides with ~(\ref{E:Up1}): $U(p_0,p_0)=U(p_0)$. Note that $U(p_0,q_0)$ is continuous in all its parameters.\\ 
To calculate the payoff of a DM who plays optimally but has an incorrect prior, simply substitute the expression for $\p1*$ from Theorem \ref{Up1} as $p'$ in Theorem \ref{UpE1}.

\subsection{Information and Payoff}\label{information}

Theorems \ref{Up1} and \ref{UpE1} express the optimal cut-off, the optimal expected payoff, and the expected payoff for a DM with incorrect prior who uses a cut-off strategy
in terms of the expected payoff of each arm,
and the quantity $\alpha$, which captures all the information-relevant parameters of the payoff processes.
In this section we justify why $\alpha$ indeed captures all the information-relevant parameters.
This justification explains why the solutions of Bolton and Harris (1995) and
Keller, Rady and Cripps (2005) have the same structure as our solution.
%
%
Suppose the DM faces a two-armed bandit problem with L\'{e}vy payoffs.
If the risky arm's type is High (resp. Low), it yields three independent L\'{e}vy processes $(X_{\high}^a (t))$, $(X_{\high}^b (t))$, $(X_{\high}^c (t))$ (resp. $(X_{\low}^a (t))$, $(X_{\low}^b (t))$, $(X_{\low}^c (t))$).

As in Section \ref{s:finite}, the notation $\mu^j_i$, $\sigma^j_i$, $\nu^j_i$, and $g^j_i$, 
represents the drift of the process $(X^j_i(t))$, the standard deviation of the Brownian motion component of $(X^j_i(t))$, the L\'{e}vy measure of the process $(X^j_i(t))$, and the expectation of the process $(X^j_i(t))$ at time $t=1$, 
respectively, where $j\in \{ a,b,c \}$, and $i\in \{ \high,\low \}$ is the arm's type.

We assume that Assumption \ref{assumption} is satisfied for the three couples ($(X_{\high}^j (t))$, $(X_{\low}^j (t))$),  $j\in \{ a,b,c \}$. 

Suppose that the DM's payoff is $X_i^a + X_i^c$, while his information is $(X_i^a, X_i^b)$. Thus, the component $X_i^a$ represents observed information of the DM that contributes to his payoff, the component $X_i^b$ represents observed information that does not contribute to the payoff, and the component $X_i^c$ represents unobserved information that contributes to the payoff. The following theorem characterizes the optimal strategy of the DM in this setup, as well as his expected payoff when he holds a possibly incorrect prior. This characterization can be used to calculate the fair price of additional information.

Let $\beta$ be the unique solution of

\begin{align}\label{E:f1}
f_{a,b} (\eta):&=
\int\vvha\left(\frac{\vvha}{\vha}\right)^{\eta}
+ \eta (\vga-\vvga)
- \vvga
+\frac{1}{2}(\eta+1)\eta \left(\frac{\aia-\bia}{\sigma^a}\right)^2\\\notag
&+\int\vvhb\left(\frac{\vvhb}{\vhb}\right)^{\eta}
+ \eta (\vgb-\vvgb)
- \vvgb
+\frac{1}{2}(\eta+1)\eta  \left(\frac{\aib-\bib}{\sigma^b}\right)^2  - r = 0\notag
\end{align} in $(0,\infty)$.

\begin{thm}\label{pricing1}
The expected payoff to a DM who holds the prior belief $q_0$, uses a cut-off strategy $\kappa '$ with cut-off point $p'$, receives and observes the payoff process $(X^a (t))$, observes but does not receive the payoff process $(X^b (t))$, and receives but does not observe the payoff process $(X^c (t))$, is as follows: if $q_0>p'$
\begin{align}\notag
V^{a,b,c}_{p'} (p_0,q_0):&=  (s-\gia-\gic)p_0\left(\frac{1-q_0}{q_0}\right)^{\beta+1} \left(\frac{p'}{1-p'}\right)^{\beta +1} \\\notag
& +(s-\ggia-\ggic)(1-p_0)\left(\frac{1-q_0}{q_0}\right)^{\beta} \left(\frac{p'}{1-p'}\right)^{\beta }\\\notag
&+p_0(\gia+\gic)+(1-p_0)(\ggia+\ggic), \notag
\end{align}
while $V^{a,b,c}_{p'} (p_0,q_0)=s$ if $q_0\leq p'$.\\
Moreover, the optimal cut-off for a DM who holds the correct prior, i.e. $q_0=p_0$, is given by
\[
 p^*= \frac{\beta (s- \ggia - \ggic)}{(\beta + 1)(\gia + \gic-s) + \beta(s - \ggia - \ggic)}.
\]
\end{thm}


Observe that if there are no processes $(X_i^b, X_i^c)$, then Theorem \ref{pricing1} reduces to Theorem \ref{UpE1}. Moreover, the parameter $\beta$ in Theorem \ref{pricing1} is the equivalent of $\alpha$ in Theorem \ref{UpE1}.

The characterization in Theorem \ref{pricing1} shows that $\beta$ incorporates all the data that is relevant to the information that the DM has. It depends on the parameters of $(X_i^a, X_i^b)$, which the DM observes, but not on the parameters of $(X_i^c)$, which the DM does not observe. Moreover, it depends on the absolute value of the difference between the drifts, $|\ai-\bi|$, the standard deviation of the Brownian motion, $\sigma$, the Radon-Nikodym derivative $\frac{\vvh}{\vh}$, and the average $\eta\vg + (1-\eta)\vvg$. These quantities help in distinguishing between the two types of the Risky arm.
Once the information-relevant parameters are summarized in $\beta$, the only relevant parameters, which affect both the optimal cut-off and the optimal expected payoff, are the expectations of the payoff-relevant processes $(g_i^a,g_i^c)$.

Theorem \ref{pricing1} can be used to find the fair price of the additional information $X^b$ and of $X^c$. That is, the price of information that does not affect the DM's payoff, as well as the price of information that affects the DM's payoff. This is done by comparing the optimal value of two bandit problems that differ in the information of the DM.

The next corollary to Theorem \ref{pricing1} states that, using cut-off strategies, additional information is more profitable.\footnote{It is well known that in one-player optimization problems, additional information cannot hurt the DM, since he can always ignore it. However, a Markovian cut-off strategy does not allow a player to forget additional information, and therefore the statement is not trivial.}

\begin{cor}\label{corollary}
Suppose that there are two decision makers, DM1 and DM2, who hold the correct prior, i.e. $q_0=p_0$, and use cut-off strategies.
DM1 receives and observes the payoff process $(X^a (t))$, observes but does not receive the payoff process $(X^b (t))$, and receives but does not observe the payoff process $(X^c (t))$.
DM2 receives and observes the payoff process $(X^a (t))$, and receives but does not observe the payoff process
$(X^c (t))$. 
Then the optimal expected payoff of DM1 is higher than the optimal expected payoff of DM2.
\end{cor}

\subsection{Optimism vs Pessimism}\label{s:optimism}

As mentioned before, the phenomenon of over confidence is common in many decision problems.
In this section we apply the result of Section \ref{incorrect prior}, and investigate who will fare better in two-armed bandit problems,
an optimist who assigns a probability higher than the true probability to the High type,
or a pessimist who assigns a probability lower than the true probability to the High type.

Suppose that there are two decision makers, DM1 and DM2,
who face independent identical copies of the decision problem.
DM1 is an \emph{optimist},
and believes that the probability of High is $p_0+\rho$, where $\rho>0$,
while DM2 is a \emph{pessimist},
and believes that the probability of High is $p_0-\rho$,
and both play optimally given their beliefs. 
If $p_0-\rho \leq \p1*$, the pessimist will always choose the safe arm, since according to his subjective belief the prior is at most the cut-off. Assume then that $p_0 -\rho >\p1*$. For every $\epsilon \in [\p1*-p_0, 1-p_0] $ denote by $V_{p_0} (\epsilon)$ the expected payoff for a DM playing optimally according to the incorrect prior $p_0 +\epsilon$, 
where $\epsilon$ may be negative.
%
It turns out that the answer regarding who will fare better, an optimist or a pessimist, depends on $\alpha$ that is defined in Eq.\ ~(\ref{eta}).\\

\begin{thm}\label{optpes}
Assume that Assumption \ref{assumption} holds.\\
1. If $\alpha>1$, then for every $\epsilon >0$ such that $\p1* < p_0 \pm \epsilon \leq 1$  we have $V_{p_0}(\epsilon) >V_{p_0}(-\epsilon) $: an optimist will fare better.\\
2. If $0 < \alpha < 1$, then for every $\epsilon >0$ such that $\p1* < p_0 \pm \epsilon \leq \frac{\alpha+2}{3}$ we have $V_{p_0}(\epsilon) > V_{p_0}(-\epsilon) $: an optimist will fare better; for every $\epsilon >0$ such that $  \frac{\alpha+2}{3} <p_0 \pm \epsilon<1$ we have $V_{p_0}(\epsilon) < V_{p_0}(-\epsilon) $: a pessimist will fare better.\\

\end{thm}
Thus, an optimist will fare better than a pessimist, unless $\alpha$ is low and $p_0-\varepsilon$ is high.
That is, a pessimist will fare better only if the following two conditions are satisfied:\\
1) The pessimist assigns high probability to the High type.\\
2) The two DM's are sufficiently patient, or it is easy to distinguish between the two types of the risky arm.\\
Otherwise, an optimist will fare better.\\\\
Since $\p1*= \frac{\alpha (s-\ggi)}{(\alpha+1)(\gi-s) + \alpha(s-\ggi)}$, the condition $\p1*<  \frac{\alpha +2}{3}$ is not always satisfied.
If $3\gi +\ggi \geq 4s$, then $\p1*<  \frac{\alpha +2}{3}$.
If $3\gi +\ggi < 4s$, then  $\p1*<  \frac{\alpha +2}{3}$ if and only if\\
$\alpha < \frac{4s-3\gi-\ggi - \sqrt{(4s-3\gi-\ggi)^2-8(\gi-\ggi)(\gi-s)}}{2(\gi-\ggi)}$ or
$\alpha > \frac{4s-3\gi-\ggi + \sqrt{(4s-3\gi-\ggi)^2-8(\gi-\ggi)(\gi-s)}}{2(\gi-\ggi)}$.\\\\

We obtain the same result if the biases of the optimist and the pessimist are geometric rather than absolute, that is, if the initial prior of the optimist is $(1+\epsilon)p_0=p_0+\epsilon p_0$ and the initial prior of the pessimist is $(1-\epsilon)p_0=p_0-\epsilon p_0$.
This is true since the absolute bias of the optimist and the pessimist is the same: $\epsilon p_0$.

\section{Future directions}

Our results call for further research. We here list few possible directions for future research.
\begin{itemize}
\item
We studied the case that the distribution of the High type dominates the distribution of the Low type in a strong sense
(see Assumptions \ref{assumption}).
These assumptions ensure that the discontinuities of the process of posterior belief
are always to one direction.
It would be interesting to know whether a similar characterization holds under only assumptions A1 and A2.
\item
We assumed that the payoff is distributed according to a L\'{e}vy process.
It would be interesting to solve the model when the payoff is distributed according to a geometric L\'{e}vy process.
\item
Bolton and Harris (1999),
Keller, Rady and Cripps (2004),
and Keller and Rady (2008)
study a strategic version of the model,
in which several decision makers face identical unknown arms.
Klein and Rady (2008) study a strategic version of the model,
in which the arms of two DMs are negatively correlated.
It will be interesting to solve these models,
when the payoff distributions are general L\'{e}vy processes.
\end{itemize}

\section{Proofs}
\subsection{Proof of Proposition \ref{convex}}
We already argued that $U$ is non-decreasing and convex. It is left to prove that it is continuous.
Since $U(p)$ is convex, it is continuous on $(0,1)$. If the DM had known the true type of the arm, his optimal payoff would have been $p\gi + (1-p)s$. Since this information is not available, $U(p)\leq p\gi + (1-p)s$. Since $U(p)\geq s$ for every $p$, this implies $\underset{p\rightarrow 0}{\lim}U(p) =s= U(0)$. Since the DM can follow the strategy that always selects the risky arm, $U(p)\geq p\gi + (1-p)\ggi$. This implies that $ \underset{p\rightarrow 1}{\lim}U(p)=\gi=U(1)$.
Thus, $U$ is continuous on $[0,1]$.
%
%
%
%
%
%
\subsection{The Control Equation}

In this section we express the optimal payoff $U(p)$ using the dynamic programming principle. This representation extends that in Bolton and Harris (1999) to the more general setup of L\'{e}vy payoff processes. we start by calculating how the belief of the DM is updated given his observations.
To simplify notation, it is convenient to divide various expressions by the standard deviation $\sigma$; set $\at = \frac{\ai}{\sigma}$, $\bt = \frac{\bi}{\sigma}$ and $\mut = \frac{\mu}{\sigma}$, and denote $\dpit :=\frac{1}{\sqrt{k}\sigma}\dpi ^k
=  \sqrt{k}\mut dt +  dZ(t),$ modulo the factor $\sqrt{\kappa}\sigma$, $\dpit$ is the contribution to the payoff of the linear drift and the Brownian motion.
Then the probability of obtaining a specific observation (that is not a jump) given the type is:
\begin{align}\notag
P(\dpit | \mut ) &=P(\dpit | \mu ) = \frac{1}{\sqrt{2\pi dt}} e^{-\frac{(\dpit (t)- \sqrt{k} \mut dt)^2}{2dt}} = C\cdot  e^{\sqrt {k}\mut \dpit - \frac{k\mut^2 (dt)^2}{2}   } \\\notag
                 &= C(1+\sqrt{k}\mut\dpit +o((dt)^2)), \notag
\end{align}
where $C$ is a constant independent of $\mu$.
Denote $\Pt ( \dpit | \mut ) : =  1+\sqrt{k}\mut\dpit$. \\\\
Let $p = P_t (\ta)$ be the belief at time $t$ that the risky arm is High. Then:
\begin{align} \notag
P_{t+dt}(\ta) &= \dfrac{P(\dpi^k |\ta) P(\dYPP|\ta)p   }{  P(\dpi^k |\ta) P(\dYPP|\ta)p +  P(\dpi^k |\tb) P(\dYPP|\tb)(1-p)}\\ \notag
        &= \dfrac{P(\dpit |\ta) P(\dYPP|\ta)p   }{  P(\dpit |\ta) P(\dYPP|\ta)p +  P(\dpit |\tb) P(\dYPP|\tb)(1-p) }\\\notag
        &= \dfrac{\Pt(\dpit |\ta) P(\dYPP|\ta)p   }{  \Pt(\dpit |\ta) P(\dYPP|\ta)p +  \Pt(\dpit |\tb)P(\dYPP|\tb)(1-p)}. \notag
\end{align}\\\\
This is the Bayesian belief updating, using the independence of the components in the L\'{e}vy-Ito decomposition, given the risky-arm type $\theta$. The next lemma expresses the change in the posterior belief over time. The lemma handles separately the case where there are no jumps in the time interval $[t,t+dt)$, and the case where there is a jump in this interval. 

\begin{lem}\label{dp}
$\\$\textbf{1.} Suppose that there are no jumps during the time interval $[t,t+dt)$. Then:
\begin{equation}\label{dP1}
dP_t:= P_{t+dt}(\ta)-P_{t}(\ta)= p(1-p) \sk (\at-\bt) \dZt - p(1-p) k (\vg-\vvg)dt,
\end{equation}

where $\dZt = \dpit - \sk (\mat)dt$ is a standard Brownian motion.\\
\textbf{2.} Suppose that during the interval $[t, t+dt)$ a jump of size $h$ occurred. Then:
\begin{equation}\label{dP2}
P_{t+dt}(\ta) = \Ph + \sk(\at-\bt)\Ph(1-\Ph)d\tilde{Z}_2,
\end{equation}
where $\Ph:= \frac{p\vh}{p\vh+(1-p)\vvh}$, and $d\tilde{Z}_2 := \dpit - \sk(\at\Ph +\bt(1-\Ph))dt$.\\
\end{lem}
The first term in the right-hand side of ~(\ref{dP2}) is the contribution of the continuous part of the payoff process to the change in the belief, while the second term is the contribution of the fact that no jump arrived. This latter contribution is negative due to Assumption A3.
If a jump of size $h$ occurs during the interval $[t,t+dt)$, then the contribution of the continuous part of the payoff process is $(\at-\bt)(\dpit - \sk(\at\Ph +\bt(1-\Ph))dt)$, and the compound Poisson process' contribution is $\Ph:= \frac{p\vh}{p\vh+(1-p)\vvh}$. The latter is the Bayesian update of the probability that the risky arm is High given that a jump of size $h$ occurred. Note that  $ p < P_h$.

\begin{proof}[\textbf{Proof of Lemma \ref{dp}}] The proof of the lemma is standard and non-inspiring. The first statement follows from a long chain of equalities. Assume that there is no jump in the interval $[t,t+dt)$
\begin{align} \notag
dP_t =& \dfrac{p(1-p) [  \Pt(\dpit |\ta) P(\dYPP|\ta) -\Pt(\dpit |\tb) P(\dYPP|\tb)   ]} {  \Pt(\dpit |\ta) P(\dYPP|\ta)p +  \Pt(\dpit |\tb) P(\dYPP|\tb)(1-p)   }\\ \notag
     =&\dfrac{p(1-p) [(1+\kt\at\dpit) ( e^{-\vg k dt} ) -   (1+\kt\bt\dpit) ( e^{-\vvg k dt} )   ]   }{  (1+\kt\at\dpit) ( e^{-\vg k dt} )p + (1+\kt\bt\dpit) ( e^{-\vvg k dt} )(1-p)    }\\ \notag
     =& \dfrac{p(1-p) [(1+\kt\at\dpit) ( 1-\vg k dt ) -   (1+\kt\bt\dpit) ( 1-\vvg k dt )   ]   }{  (1+\kt\at\dpit) ( 1-\vg k dt )p + (1+\kt\bt\dpit) ( 1-\vvg k dt )(1-p)    } \\ \notag
     =& \dfrac{ p(1-p)[ \sk(\at-\bt)\dpit -(\vg-\vvg)kdt  ]    }{ 1+ \sk (\mat)\dpit -k(\mvt)dt     } \\ \notag
     =& p(1-p) [ \sk(\at-\bt)\dpit -(\vg-\vvg)kdt  ]\cdot \\ \notag
     &  \cdot [1- \sk(\mat)\dpit + k(\mvt)dt + k(\mvt)^2 dt    ] \\ \notag
     =& p(1-p)[ \sk(\at -\bt) (\dpit -\sk (\mat)dt) - k(\vg -\vvg)dt    ] \\ \notag
     =& p(1-p) \sk (\at-\bt) \dZt - p(1-p) k (\vg-\vvg)dt. \\ \notag
\end{align}
In  the calculations we used the fact that $\dZt = \dpit - \sk (\mat)dt$ is a standard Brownian motion (see Bolton and Harris (1999)), and the Brownian motion properties:
$dZ^2 = dt$,  and $dZ dt = 0 $. We also ignored terms of order $(dt) ^ {3/2}$ and above.\\\\\\\\
We now prove the second statement.
\begin{align}\notag
P_{t+dt}(\ta) &=  \dfrac{p\vh ( 1+\sk \at\dpit)     }{ [ \vh p +\vvh (1-p) ] + \sk[\vh\at p + \vvh \bt (1-p)  ]\dpit     }\\ \notag
         &= \dfrac{ \Ph (1+\sk\at\dpit)  }{1+\sk[\at\Ph + \bt(1-\Ph) ]\dpit} \\ \notag
         &= \Ph (1+\sk\at\dpit)[ 1-\sk(\Ph \at +(1-\Ph)\bt)\dpit +k(\Ph \at +(1-\Ph)\bt)^2  dt  ] \\ \notag
         &= \Ph [1+\sk(\at-\bt)(1-\Ph)(\dpit - \sk(\at\Ph +\bt(1-\Ph))dt   )] \\ \notag
         &= \Ph + \sk(\at-\bt)\Ph(1-\Ph)d\tilde{Z}_2.  \notag
\end{align}
         %
         %
         %
         %
\end{proof}
%
%
%
%
%
%
%
%
%
%
%
%
%
%
%
%
%
%
We now formulate the control equation that describes the optimal payoff.\\\\
\begin{equation}
U(p)=\underset{k\in[0,1]}{max} \{[(1-k)s+ k(p(\vg \Hi + \ai) +(1-p)(\vvg \HHi +\bi))]rdt + e^{-rdt} E[U(p+dp)]\}, \tag {CE}
\end{equation}
$\\$where $k$ is the control variable. The first term within the maximization is the expected instantaneous payoff, and the second term is the discounted expected continuation payoff.
The following lemma provides a more convenient form to the control equation, in terms of the derivatives of $U$.

\begin{lem}\label{LCE} The following equality holds:
\begin{align} \tag{CE2}
 U(p) =  \underset {k\in[0,1]}{max} &\left\{ (1-k)s+ k(p(\vg \Hi + \ai) +(1-p)(\vvg \HHi +\bi))\right.  \\ \notag
 &+\frac{1}{r} \left[ k \int  (\mvh)U\left(\frac{p\vh}{\mvh}\right)  \right. \\ \notag
& - kp(1-p)(\vg-\vvg)U'(p) - k(\mvg) U(p)\\\notag
&+ \left. \left.\frac{1}{2}  k U''(p)p^2(1-p)^2 (\at-\bt)^2 \right] \right\}, \quad p-a.s.
\end{align}
\end{lem}

\begin{proof}[\textbf{Proof of Lemma \ref{LCE}}]

Since $U(p)$ is convex, $U(p)$ is twice differentiable $p$-a.s. (in the sense of the Lebesgue measure).
With probability $[p(1-k\vg dt) + (1-p)(1-k\vvg dt)]$ there are no jumps in the interval $[t,t+dt)$. In this case, using the Taylor expansion of $U$, and ignoring terms of order $(dt)^{3/2}$ and higher, we obtain that the optimal payoff is $U(p+\check{d}p)=U(p) + U'(p) \check{d}p +\frac{1}{2} U''(p) (\check{d}p)^2$ a.s.,\footnote{Since we can ignore terms or order $(dt)^{3/2}$ and higher, it is sufficient to consider the Taylor expansion up to the second derivative.} where $\check{d}p$ is given by the right-hand side of ~(\ref{dP1}).\\\\
With probability  $[p k\vh dt + (1-p)k\vvh dt]$ there is a jump of size $h$, and the optimal payoff is $U(\Ph + \sk(\at-\bt)\Ph(1-\Ph)d\tilde{Z}_2) = U(\Ph)+ U'(\Ph)\hat{d}\Ph + \frac{1}{2} U''(\Ph)\hat{d}\Ph$,
where $\hat{d}\Ph: =\sk(\at-\bt)\Ph(1-\Ph)d\tilde{Z}_2$.  \\
$\\$During the subsequent calculations we use the following Eqs.\ ~(\ref{Edp}), ~(\ref{Edp2}), ~(\ref{EdPh}) and ~(\ref{EdPh2}), that can be derived from Lemma \ref{dp}:
\begin{equation}\label{Edp}
E[\check{d}p] = -kp(1-p)(\vg-\vvg)dt.
\end{equation}
This is the expectation of the change in the belief, given that no jump occurred during the interval $[t,t+dt)$.
\begin{equation}\label{Edp2}
E[\check{d}p^2] = kp^2 (1-p)^2 (\at-\bt)^2 dt.
\end{equation}
This is the second moment of the change in the belief, given that no jump occurred during the interval $[t,t+dt)$.
\begin{equation}\label{EdPh}
E[\hat{d}\Ph] = C_1 \cdot dt.
\end{equation}
 This is the expected contribution of the Brownian motion part to the posterior belief, given that a jump of size $h$ occurred during the interval $[t,t+dt)$.
 \begin{equation}\label{EdPh2}
 E[\hat{d}\Ph^2] = C_2 \cdot dt.
 \end{equation}
 This is the second moment of the contribution of the Brownian motion part to the posterior belief, given that a jump of size $h$ occurred during the interval $[t,t+dt)$.
In Eqs.\ ~(\ref{EdPh}) and ~(\ref{EdPh2}), $C_1$ and $ C_2$ are constants.
Using the above notation, we obtain from (CE):
\begin{align}\label{ce1}
U(p) = \underset{k\in[0,1]}{max} &\left\{  [(1-k)s+ k(p(\vg \Hi + \ai) +(1-p)(\vvg \HHi +\bi))]rdt \right. \\ \notag
&+ (1-rdt)\left[ kdt \left[ \int  \left[U\left(\frac{p\vh}{\mvh}\right)\right.\right.\right. \\\notag
&\left.\left.\left.+ U'\left(\frac{p\vh}{\mvh}\right) C_1 dt \right.\right.\right.\\ \notag
&+   \left.\left.U''\left(\frac{p\vh}{\mvh}\right) C_2 dt\right](\mvh)  \right]  \\ \notag
&+ \left.\left.\left[1-k(\mvg) dt\right]\left[U(p) + U'(p)E[\check{d}p]+\frac{1}{2}U''(p)E[\check{d}p^2]\right] \right]\right\}, \quad p-a.s.\notag
\end{align}
The second, third and fourth lines in ~(\ref{ce1}) represent the expected continuation payoff given a jump of size $h$ occurred during the time interval $[t,t+dt)$, and the fifth line represents the expected continuation payoff given no jump occurred during the time interval $[t,t+dt)$.
Using $(dt)^2=0$, several of the terms in ~(\ref{ce1}) vanish, and we obtain:
\begin{align}\label{ce2}
U(p) =  \underset{k\in[0,1]}{max} & \left\{ [(1-k)s+ k(p(\vg \Hi + \ai) +(1-p)(\vvg \HHi +\bi))]rdt\right.  \\ \notag
&+ kdt \int U\left(\frac{p\vh}{\mvh}\right) (\mvh)  \\ \notag
 &+ U(p) - kdt p(1-p) (\vg-\vvg)U'(p)dt - k(\mvg) U(p)dt -rU(p)dt \\ \notag
&+ \left.\frac{1}{2} k U''(p)p^2(1-p)^2 (\at-\bt)^2 dt  \right\}, \quad p-a.s.\\ \notag
 \end{align}
Eliminating $U(p) $ from both sides, and dividing by $dt$, we obtain (CE2) after simple algebraic manipulations, as desired.
\end{proof}

From Eq.\ ~(\ref{ce1}) it follows that the contribution of the continuation payoff given that a jump of size $h$ occurred during the time interval $[t,t+dt)$ is
\begin{align}\label{ce3}
 kdt &\left[ \int  \left[U\left(\frac{p\vh}{\mvh}\right)\right.\right.  \left.\left.\left.+ U'\left(\frac{p\vh}{\mvh}\right) C_1 dt \right.\right.\right.\\ \notag
&+   \left.\left.U''\left(\frac{p\vh}{\mvh}\right) C_2 dt\right](\mvh)  \right].\notag
\end{align}
The parameters of the Brownian motion affect ~(\ref{ce3}) only through $C_1$ and $C_2$, and since $C_1$ and $C_2$ do not appear in Eq.\ ~(\ref{ce2}), it follows that if a jump occurs during the time interval $[t,t+dt) $, the information from the compound Poisson process has more impact than the information of the Brownian motion.\\

$\\$According to (CE2), the payoff is the maximum over the control variable $k$ of the expectation of the current flow payoff  $\left[(1-k)s+ k(p(\vg \Hi + \ai) +(1-p)(\vvg \HHi +\bi))\right]$ plus the discounted value of the continuation payoff
\begin{align}\notag
\frac{1}{r} &\left[ k \int  (\mvh)U\left(\frac{p\vh}{\mvh}\right)  \right. \\ \notag
                   & - kp(1-p)(\vg-\vvg)U'(p)+  \left.\frac{1}{2}  k U''(p)p^2(1-p)^2 (\at-\bt)^2 - k(\mvg) U(p) \right].\notag
\end{align}
A solution $\kappa$ to this maximization problem must satisfy:\\\\
\begin{equation}\label{k}
k =
\begin{cases}
0               &\text{if $b(p,U)<s - [p\gi +(1-p)\ggi]$}, \\
\in[0,1]        &\text{if $b(p,U)=s - [p\gi +(1-p)\ggi]$},\\
1               &\text{if $b(p,U)>s - [p\gi +(1-p)\ggi]$},
\end{cases}
\end{equation}
$\newline$
where
\begin{align}\notag
b(p,U) = \frac{1}{r} &\left[  \int  (\mvh)U\left(\frac{p\vh}{\mvh}\right) \right.\\\notag
                     &- p(1-p)(\vg-\vvg)U'(p) - (\mvg) U(p)\\\notag
                     &\left.+\frac{1}{2}   U''(p)p^2(1-p)^2 (\at-\bt)^2 \right].\notag
\end{align}
The function within the maximization in (CE2) is linear in $k$. Therefore, it achieves its maximum at $k=1$ or $k=0$, $p$-a.s. From Proposition \ref{convex}, $U(p)$ is non-decreasing and continuous, and therefore there is $\p1* $ such that $U(p)=s$ for every $p\leq \p1*$. Thus $k=0$ is optimal for $p<\p1*$. For every $p>\p1*$, we have $U(p)>s$, so that in this case $k=1$ is optimal $p$-a.s.\footnote{Recall that Eq.\ (CE2) is satisfied $p$-a.s., since $U'(p)$ and $U''(p)$ exist $p$-a.s. In the next Section we show that the optimal strategy is in fact a cut-off strategy.} \\
\subsection{Characterizing the optimal strategy and the value}\label{characterizing}
%
When it is optimal to play safe, that is, when the optimal solution of (CE) is $k^*=0$, we have  $U(p)=s$.
When it is optimal to play risky, that is, when the optimal solution of (CE) is $k^*=1$,  it follows from Lemma \ref{LCE} that $U(p)$ solves the following functional differential equation:
\begin{align} \tag{FDE}
U(p) &= p(\vg \Hi + \ai) +(1-p)(\vvg \HHi +\bi)  \\ \notag
      &+\frac{1}{r} \left[  \int  (\mvh)U\left(\frac{p\vh}{\mvh}\right) \right.\\ \notag
      & - p(1-p)(\vg-\vvg)U'(p)- (\mvg) U(p)\\\notag
      &+\left.\frac{1}{2}   U''(p)p^2(1-p)^2 (\at-\bt)^2 \right],\; a.s.\;\, in\;\, (\p1*,1).
\end{align}
A solution $U(p)$ for this equation must be smooth (Friedman (1969), p.56).\footnote{ To see that the conditions of Friedman (1969) are satisfied, substitute $f(p)=\int  (\mvh)U\left(\frac{p\vh}{\mvh}\right) - U(p)(r+(\mvg))$. Since $U(p)$ is continuous, so is $f$, and we get $U\in C^2$, as Friedman (1969) requires.}
Therefore, $U(p) $ satisfies Eq.\ (FDE) in $(\p1*,1)$ always, and $k=1$ is optimal in $(\p1*,1)$. Hence, there is an optimal cut-off strategy $\kappa^* $ with cut-off point $\p1*$.\\

The next lemma suggests one solution to Eq.\ (FDE). To prove it, substitute the expression for $U(p)$ defined in Eq.\ ~(\ref{expressionUp}) below into Eq.\ (FDE).
Recall that $$f(\eta) =  \int\vvh\left(\frac{\vvh}{\vh}\right)^{\eta} + \eta (\vg-\vvg)-\vvg +\frac{1}{2}(\eta+1)\eta (\at-\bt)^2 - r $$ (see ~(\ref{eta})), that $\alpha$ is the unique solution of the equation $f(\eta)=0$ in $(0,\infty)$, and that $\p1*$ and $C_\alpha$ were defined in the statement of Theorem \ref{Up1}.

\begin{lem}\label{solution}
One smooth solution to Eq.\ (FDE) is
\begin{equation}\label{expressionUp}
U(p) =  p\gi+(1-p)\ggi +  C_\alpha(1-p)\left( \frac{1-p}{p}\right)^{\alpha},
\end{equation}
where $\;\alpha\in (0,\infty)$ 
solves the equation $f(\eta) = 0 $.
\end{lem}
In order to see that the function $U(p)$, defined above, actually solves (FDE), we use the fact that $P_h=\frac{p\vh}{\mvh}\geq \p1*$ for every $p>\p1*$, $\nu_1 -$a.s., which is equivalent to Assumption A3. Thus, the form of our solution crucially depends on this assumption.

In fact, one can verify that $$U(p) =   p\gi +(1-p)\ggi +  C(1-p)\left( \frac{1-p}{p}\right)^{\alpha} +  D(1-p)\left( \frac{1-p}{p}\right)^{\beta}$$ solves Eq.\ (FDE), where $\alpha$ is as in the statement of Lemma \ref{solution}, and $\beta$ is the unique solution\footnote{In fact, such a solution $\beta $ must be smaller then $-1$.} of $f(\eta) = 0$ in $(-\infty,0)$.
%
%
%
The following lemma assures that $\alpha$ is well defined.
\begin{lem}\label{alpha}
The equation $f(\eta) =0$ admits a unique solution in the interval $(0,\infty)$.\\
\end{lem}
\begin{proof}[\textbf{Proof}]
$\newline$
The function $f$ is a continuous function that satisfies $f(0)<0$ and $f(\infty) = \infty$. To show that $f(\eta)=0$ has a unique solution, it is therefore sufficient to prove that $f$ is increasing in $\eta$. Note that if $\ai\neq\bi$, then
$\frac{1}{2}(\eta+1)\eta (\at-\bt)^2 - r -\vvg$ is increasing in $\eta$, and constant otherwise. It remains to prove that if $\vi\neq\vvi$ (i.e. $\vi (\mathbb{R}\backslash\{0\})>\vvi(\mathbb{R}\backslash\{0\})$) , then $\int \vvh\left(\frac{\vvh}{\vh}\right)^{\eta} + \eta (\vg-\vvg)$ is increasing in $\eta$. Since
$$\int \vvh\left(\frac{\vvh}{\vh}\right)^{\eta} + \eta (\vg-\vvg) =
\int\left[\vvh\left(\frac{\vvh}{\vh}\right)^{\eta} + \eta (\vh-\vvh) \right],$$\\
and $$\underset{ \{ h|\frac{\vvh}{\vh}=1  \}   }{\int}\left[\vvh\left(\frac{\vvh}{\vh}\right)^{\eta} + \eta (\vh-\vvh) \right]=0,$$
it is sufficient to prove that for $\vi-$a.e. $h\in \{ h| \frac{\vvh}{\vh}<1  \} $,
$$g_h(\eta) = \vvh\left(\frac{\vvh}{\vh}\right)^{\eta} + \eta (\vh-\vvh) $$ is increasing in $\eta$. Now,
\begin{align}\notag
g_h'(\eta) &= -\vvh \ln\left(\frac{\vh}{\vvh} \right)\left(\frac{\vvh}{\vh}\right)^{\eta} +  (\vh-\vvh)  \\\notag
&> -\vvh \ln\left(\frac{\vh}{\vvh} \right) +  (\vh-\vvh) > 0,
 \end{align}
where the first inequality holds since $\eta>0$ and by Assumption A3. The second inequality holds since $-\ln(x) +x-1> 0$ for every $x\neq 1$. Therefore, $g_h(\eta)$ is increasing, as desired.
\end{proof}
%
%
%
%
%
%
%
%
%
We now prove that Eq.\ (FDE) has a unique solution.

\begin{lem}\label{uniqueness}
For every $ p_1 < p_2 $, and every  $  u_1,u_2\in \mathbb{R} $, there is a unique solution $U(p)$ satisfying Eq.\ (FDE) in the interval $(p_1,p_2)$ with the boundary conditions $U(p_1)=u_1,\; U(p_2)=u_2$.
\end{lem}

\begin{proof}[\textbf{Proof}]
$\newline$
Since Eq.\ (FDE) is a non-homogenous linear equation in $U$, if there are two solutions of Eq.\ (FDE), then their difference is a solution of the homogenous version of Eq.\ (FDE). To prove the lemma, it is therefore sufficient to fix a solution $W$ of the homogenous version of Eq.\ (FDE) that satisfies $W(p_1) = W(p_2)=0$ and to prove that $W\equiv 0$.

Suppose that $W$ achieves its maximum at $\hat{p}$. Then $W'(\hat{p})=0, \,$ therefore: \\
\begin{align}\notag
W(\hat{p}) &= \frac{1}{r} \left[  \int  (\mvhgag)W\left(\frac{\hat{p}\vh}{\mvhgag}\right) \right.\\ \notag
      &+\left.\frac{1}{2}   W''(\hat{p})\hat{p}^2(1-\hat{p})^2 (\at-\bt)^2 - (\mvggag) W(\hat{p})\right] .\\\notag
\end{align}

Moreover, since the maximum is achieved at $\hat{p}$, $ W''(\hat{p})\leq 0$, simple algebraic manipulations imply that:
\begin{align}\notag
(r+ \mvggag)W(\hat{p}) &=  \int  (\mvhgag)W\left(\frac{\hat{p}\vh}{\mvhgag}\right) \\\notag
                        &+ \frac{1}{2}W''(\hat{p})\hat{p}^2(1-\hat{p})^2 (\at-\bt)^2 \\\notag
                        &\leq W(\hat{p}) \int  (\mvhgag)\\\notag
                        &= (\mvggag)W(\hat{p}).\\\notag
\end{align}
Since $r> 0 $ we conclude that $W(\hat{p}) = 0$. A similar argument shows that the minimum of $W$ in $(p_1,p_2)$ is $0$, so that $W(p)\equiv 0 $ on $(p_1,p_2)$, as desired.
\end{proof}
%
%
%
%
%
As mentioned before, there is an optimal cut-off strategy with corresponding payoff $U$. Lemmas \ref{solution}, \ref{alpha} and \ref{uniqueness} prove that $U$ is the unique solution of Eq.\ (FDE). We now prove Theorem \ref{Up1}, which provides an explicit form to the optimal strategy and to the payoff function.

\begin{proof}[\textbf{Proof of Theorem \ref{Up1}}]
$\newline$
Recall that $\kappa^*$ is the optimal cut-off strategy (with cut-off point $\p1*$). To complete the proof of the theorem, we provide an explicit expression to $\p1*$ and to $U$. To this end, we first prove that the right-hand derivative of $U$ at $\p1*$ is 0.

%
Let $U'_R (\p1*)=\underset{\varepsilon\rightarrow 0^+}{\lim}\frac{U(\p1*+\varepsilon)-U(\p1*)}{\varepsilon} $ be the right derivative of $U$ at the cut-off point $\p1*$. Since $U$ is convex, $U'_R(\p1*)$ is well defined. Since 
$U$ is non-decreasing, $U'_R (\p1*)\geq 0$. We now prove that  $U'_R (\p1*)\leq 0$.
For every $q_0\in[0,1]$, let $\kappa (q_0)$ be the strategy that plays as $\kappa$, assuming the prior belief is $q_0$ rather then $p_0$. Define $M_{\kappa} := \int_0^{\infty} r e^{-rt}dY^{\kappa} (t)$, the discounted payoff under the strategy $\kappa$. 
\[
 \underset{\epsilon\rightarrow 0^+ }{\lim}E\left[\left.M_{\kappa^* (\p1*+\epsilon)}\right|\theta\right] =s\quad \forall \theta\in\{\ta, \tb\}.
 \]
%
Indeed, as explained before, the posterior belief will drop below $p_0$ in an infinitesimal time interval around $0$. Therefore,\footnote{By the same argument we get that $V_{\kappa '}(p)$  is continuous in $p'$, where $\kappa '$ is a cut-off strategy with cut-off point $p'$. } as $\epsilon$ goes to $0$, the probability that the DM stops ``quite fast" under $\kappa^* (\p1*+\epsilon)$ goes to $1$.\\\\

Since $\kappa^*$ is the optimal cut-off strategy, and since it is independent of the prior belief $p_0$, we deduce that for every $p\in[0,1] $, $U(p)=V_{\kappa^* (p)} (p)$. Therefore:

\begin{align}\label{right1}
U'_R (\p1*)  &=
\underset{\epsilon\rightarrow 0^+}{\lim} \frac{U(\p1*+\epsilon)-U(\p1*)}{\epsilon}\\\notag
&=\underset{\epsilon\rightarrow 0^+}{\lim} \frac{V_{\kappa^* (\p1*+\epsilon)}(\p1*+\epsilon)-V_{\kappa^* (\p1*)}(\p1*)}{\epsilon}\\\notag
&=\underset{\epsilon\rightarrow 0^+}{\lim}\dfrac{1}{\epsilon}\left[ (\p1*+\epsilon) E[M_{\kappa^* (\p1*+\epsilon)}|\ta]\right.
\left.+(1- \p1*-\epsilon) E[M_{\kappa^* (\p1*+\epsilon)}|\tb]\right]\\\notag
&\qquad\quad \left.-\p1* E[M_{\kappa^* (\p1*)}|\ta]
-(1- \p1*) E[M_{\kappa^* (\p1*)}|\tb]\right]\\\notag
&=\underset{\epsilon\rightarrow 0^+}{\lim}\left\{\dfrac{1}{\epsilon}\left[ \p1*  E[M_{\kappa^* (\p1*+\epsilon)}|\ta]\right.\right.
+(1-\p1*)  E[M_{\kappa^* (\p1*+\epsilon)}|\tb]\\\notag
&\qquad\quad \left.-\p1* E[M_{\kappa* (\p1*)}|\ta]
-(1-\p1*) E[M_{\kappa^k (\p1*)}|\tb]\right]\\\notag
&\qquad\quad\left. +E[M_{\kappa^* (\p1*+\epsilon)}|\ta] - E[M_{\kappa^* (\p1*+\epsilon)}|\tb]\right\}\notag.
\end{align}
By the optimality of $U(p)$, $U(\p1*)\geq V_{\kappa^*(\p1*+\epsilon)}(\p1*)$, and therefore,\\
\begin{align}\label{right2}
&\underset{\epsilon\rightarrow 0^+}{\lim}\dfrac{1}{\epsilon}\left[ \p1*  E[M_{\kappa^* (\p1*+\epsilon)}|\ta]
+(1-\p1*)  E[M_{\kappa^* (\p1*+\epsilon)}|\tb]
-\p1* E[M_{\kappa^* (\p1*)}|\ta]
-(1-\p1*) E[M_{\kappa^* (\p1*)}|\tb]\right] \\\notag
&= \underset{\epsilon\rightarrow 0^+}{\lim}\dfrac{1}{\epsilon}
[V_{\kappa^*(\p1*+\epsilon)}(\p1*)- V_{\kappa^*(\p1*)}(\p1*)]=
\underset{\epsilon\rightarrow 0^+}{\lim}\dfrac{1}{\epsilon}
[V_{\kappa^*(\p1*+\epsilon)}(\p1*)- U(\p1*)]
\leq 0,
\end{align} and
\begin{equation}\label{right3}
\underset{\epsilon\rightarrow 0^+}{\lim} \left(E[M_{\kappa^* (\p1*+\epsilon)}|\ta] - E[M_{\kappa^* (\p1*+\epsilon)}|\tb]\right) = 0.
\end{equation}
Substituting ~(\ref{right2}) and ~(\ref{right3}) in ~(\ref{right1}) we deduce that $U'_R (\p1*) \leq 0$.\\

As mentioned in Lemma \ref{convex}, $U(\p1*)=s,$ and $  U(1)=\gi$. By Lemma \ref{uniqueness}, the unique solution of Eq.\ (FDE) on $(\p1*,1]$ is $U(p) =  p\gi+(1-p)\ggi +  C_{\alpha}(1-p) \left( \frac{1-p}{p}\right)^{\alpha} $.\\
By imposing $U(\p1*)=s$ and $ U'_R(\p1*)=0$ we get $\p1*= \frac{\alpha (s-\ggi)}{(\alpha+1)(\gi-s) + \alpha(s-\ggi)},$ and $C_{\alpha}= \frac{s-\ggi - \p1* (\gi -\ggi)}{(1-\p1*) ( \frac{1-\p1*}{\p1*})^{\alpha}}$. Uniqueness follows by ~(\ref{k}).\\
\end{proof}
%
%
%
%
%
%
%
%
\subsection{Incorrect Prior}\label{incorrect}
To find the expected discounted payoff for a DM who plays according to an incorrect prior $q_0$, we present here a condition which is equivalent to ~(\ref{posteriorq}).
Substituting $\pit = \mu t +\sigma Z(t)$ and  $\mut = \frac{\mu}{\sigma}$ in ~(\ref{posteriorq}), we get:\\
$Z(t) + \left[ \frac{2\mut-\at-\bt}{2} - \frac{\vg-\vvg}{\at-\bt} \right]t >  -\frac{1}{\at-\bt} \left[\ln\left(\frac{q_0}{1-q_0}\right)-\ln\left(\frac{p'}{1-p'}\right)\right]-\frac{1}{\at-\bt} \sum_{t-} \lnvhvv.$
%
It follows that the DM selects the risky arm until the first time $t$ is satisfied:
\begin{align}\label{Bmut}
B^{\mu}(t) &:= Z(t) + \left[ \frac{2\mut-\at-\bt}{2} - \frac{\vg-\vvg}{\at-\bt} \right]t \\\notag
&\leq  -\frac{1}{\at-\bt} \left[\ln\left(\frac{q_0}{1-q_0}\right)-\ln\left(\frac{p'}{1-p'}\right)\right]-\frac{1}{\at-\bt} \sum_{t-} \lnvhvv .\notag
\end{align}\\

\begin{proof}[\textbf{Proof of Theorem \ref{UpE1}}]

If the prior belief of the DM, $q_0$, satisfies $q_0 \geq p'$ then there is a bijection relation between $E$ and $q_0$, 
If $q_0 \leq p'$, then the DM always chooses the safe arm, which is equivalent to $E=0$. Therefore, we will use the notation $U(p_0,E)$ instead of $U(p_0,q_0)$ when the former is more convenient. We now prove that under Assumption \ref{assumption}, for every $p_0\in [0,1]$ and every $ E\in[0,\infty)$, the payoff of the DM is\\
$U(p_0,E) =  p_0\gi +  (1-p_0)\ggi + (s-\gi)p_0 e^{-(\at-\bt)(\alpha+1) E}+ (s-\ggi)(1-p_0)e^{-(\at-\bt)\alpha E}. $\\

Using Eq.\ ~(\ref{Bmut}), we construct an integral equation, to find the utility function for a DM who has a prior belief $q_0$. 
Let $\tau$ be the stopping time of the first jump. Let $T$ be the first time $t$ satisfying ~(\ref{Bmut}).
The DM chooses the risky arm until the stopping time $T$, and then he switches to the safe arm.
The calculations use dynamic programming in which the continuation payoff is determined by the time of the first jump, $\tau$, and the value of the continuous part of the payoff at that time.\\\\

\textbf{Notation and Formulas.} 
The proof requires computations that rely on some results on Brownian motion. In this subsection we provide these results, that are derived using Borodin and Salminen (1996) p.197 - 223. Recall that, $\tilde{\mu}=\frac{\mu}{\sigma}$ is determined by $\theta$.
Notice that $B^{\mu} (t)$ is a standard Brownian motion with drift $F_{\mu} :=\left[\dfrac{2\mut-\at-\bt}{2} - \dfrac{\vg-\vvg}{\at-\bt}\right]t$.
Define $\fmu :=  P_\theta (\Bmutau \in dx, \tau =t  | \tau<T)$. This is the probability that the first jump occurs in the interval $[t,t+dt)$, and $B^{\mu}$ belongs to $[x,x+dx)$, given a jump occurs before the DM switches to the safe arm.

Denote
%
\begin{align}\label{1}
&p_{t,h,x} := P(\ta | \tau<T, \tau =t, \Bmutau \in dx , h) \\\notag
       &= \frac{ P( \tau<T, \tau =t, \Bmutau \in dx, h |\ta )P(\ta)}{P( \tau<T, \tau =t, \Bmutau \in dx, h |\ta )P(\ta) + P( \tau<T, \tau =t, \Bmutau \in dx, h |\tb )P(\tb)}  \\\notag
       &=\frac{p_0 \Patau \fa \vh/\vg}{p_0 \Patau \fa \vh/\vg + (1-p_0) \Pbtau \fb \vvh/\vvg}.        \notag
\end{align}
This is the posterior belief that the type is $\theta_1$, given that (a) the first jump that occurred in the time interval $[t,t+dt)$ has size $h$; (b) it occurred before the DM switched to the safe arm; and (c) the Brownian motion with drift, $B^{\mu}(t) $, is in the interval $[x,x+dx)$.\\\\
The probability that the DM switches to the safe arm before the first jump appeared is
\begin{align}\label{2}
P_{\theta}(\tau>T) = P_{\theta}(\underset{0<s<\tau}{\inf}\Bsmu \leq -E )=e^{-E (\Fmu + \sqrt{2\bar{\nu}+\Fmu^2})}.\\\notag
\end{align}
The expected discounted payoff from the continuous part of the risky arm, until the switching time to the safe arm, given the DM switched before the first jump, is
\begin{align}\label{3}
E_{\theta}\left[ \left.\int_0^{T}re^{-rt} d\pita \right| \tau > T \right] &= \mu E_{\theta}\left[ \left.\int_0^{T}re^{-rt} dt \right| \tau > T \right] + \sigma E_{\theta}\left[ \left.\int_0^{T}re^{-rt} dZ_t \right| \tau > T \right] \\\notag
                             &=\mu E_{\theta}[1-e^{-rT} |\tau > T] = \mu (1-E_{\theta}[e^{-rT}| \tau > T]).\\\notag
\end{align}
The expected discounted payoff from the safe arm, after the switching time to the safe arm, given the DM switched before the first jump, is
\begin{align}\label{4}
E_{\theta}\left[\left. \int_T^{\infty} re^{-rt} sdt \right| \tau > T\right] = sE_{\theta}\left[ e^{-rT}| \tau > T  \right].\\\notag
\end{align}
The expected discounted payoff from the continuous part of the risky arm, until the first jump occurs, given the first jump occurred before the switching time, is
\begin{align}\label{5}
E_{\theta}\left[ \left.\int_0^{\tau}re^{-rt} d\pita \right| \tau < T \right] = \mu (1-E_{\theta}\left[e^{-r\tau}|\tau<T\right]).\\\notag
\end{align}
The expression on the right-hand side of ~(\ref{4}), can be re-written using the following list of equalities:
\begin{align}\label{6}
E_{\theta}[e^{-rT}| \tau>T] &= \frac{E_{\theta}[e^{-rT}, \tau>T]}{P_{\theta}(\tau>T)} = \frac{\int e^{-rt_1} P_{\theta}(T\in dt_1 , t_1<\tau)dt_1  }{P_{\theta}(\tau>T)}  \\\notag &= \frac{\int e^{-rt_1} P_{\theta}(T\in dt_1)P_{\theta}( t_1<\tau)dt_1  }{P_{\theta}(\tau>T)}
=\frac{\int e^{-rt_1} P_{\theta}(T\in dt_1)e^{-\bar{\nu} t_1}dt_1  }{P_{\theta}(\tau>T)}\\\notag
&=\frac{\int e^{-(r + \bar{\nu} )t_1} P_{\theta}(T\in dt_1) dt_1 }{P_{\theta}(\tau>T)} =\frac{ E_{\theta}[e^{-(r+\bar{\nu})T}]}{P_{\theta}(\tau>T)} \\\notag
&= \frac{e^{-E(\Fmu + \sqrt{2(r+\bar{\nu})+\Fmu^2} )}}{P_{\theta}(\tau>T)} =\frac{P_{\theta}(\tau^r >T)}{P_{\theta}(\tau>T)},\\\notag
\end{align}
where $P_{\theta}(\tau^r >t) = e^{-(r+\bar{\nu})t}$, and
$P_{\theta}(\tau^r>T) =e^{-E (\Fmu + \sqrt{2(\bar{\nu}+r)+\Fmu^2})}$. Similarly, the expression on the right-hand side of ~(\ref{5}) can be re-written as follows:
\begin{align}\label{7}
E_{\theta}[e^{-r\tau}| \tau<T] &= \int e^{-rt_1} P_{\theta}(\tau \in dt_1 |\tau<T) dt_1
=\int \frac{e^{-rt_1} P_{\theta}(\tau \in dt_1 ,t_1<T)}{P_{\theta}(\tau<T)} dt_1 \\\notag
&= \int \frac{e^{-rt_1} P_{\theta}(\tau \in dt_1)P_{\theta}(t_1<T)}{P_{\theta}(\tau<T)} dt_1
=\int \frac{e^{-rt_1} \vb e^{-\vb t_1} P_{\theta}(t_1<T)}{P_{\theta}(\tau<T)} dt_1\\\notag
&=\frac{\vb}{\vb + r}  \int (\vb + r)\frac{e^{-(r+\vb)t_1} P_{\theta}(t_1<T)}{P_{\theta}(\tau<T)} dt_1 \\\notag
&=\frac{\vb}{\vb + r}  \int \frac{P_{\theta}(\tau^r \in dt_1) P_{\theta}(t_1<T)}{P_{\theta}(\tau<T)} dt_1 =\frac{\vb}{\vb + r} \frac{P_{\theta}(\tau^r <T)}{P_{\theta}(\tau<T)}.
\end{align}
In the proof of Theorem \ref{UpE1} we use the following two identities:
\begin{align}\label{8}
\int_0^{\infty}& \fmu e^{-\gamma t} dt = \frac{\int e^{-\gamma t} P_{\theta}
\left(\underset{0<s<t}{\inf}\Bsmu \geq -E , \Btmu \in dx ,\tau \in dt\right) }{P_{\theta}(\tau<T)} \\\notag
&=\frac{\int e^{-\gamma t} P_{\theta}\left(\underset{0<s<t}{\inf}\Bsmu \geq -E , \Btmu \in dx \right)P_{\theta}(\tau \in dt )dt}
{P_{\theta}(\tau<T)} \\\notag
 &=\frac{\int e^{-\gamma t} P_{\theta}\left(\underset{0<s<t}{\inf}\Bsmu \geq -E , \Btmu \in dx\right) \vb e^{-\vb t}dt }{P_{\theta}(\tau<T)} \\\notag
&= \frac{\vb}{\vb + \gamma}\frac{\int  (\vb+\gamma )e^{-(\gamma+\vb) t} P_{\theta}\left(\underset{0<s<t}{\inf}\Bsmu \geq -E , \Btmu \in dx\right)dt}{P_{\theta}(\tau<T)} \\\notag
 &=
\frac{\vb}{\vb + \gamma}\frac{\int   P_{\theta}\left(\underset{0<s<t}{\inf}\Bsmu \geq -E , \Btmu \in dx\right)P_{\theta}(\tau^{\gamma}\in dt)dt}{P_{\theta}(\tau<T)} \\\notag
&=\frac{\vb}{\vb + \gamma}\frac{\int   P_{\theta}\left(\underset{0<s<t}{\inf}\Bsmu \geq- E , \Btmu \in dx , \tau^{\gamma}\in dt\right)dt}{P_{\theta}(\tau<T)} \\\notag
 &=
\frac{\vb}{\vb + \gamma}\frac{  P_{\theta}( B^{\mu}(\tau^{\gamma})\in dx, \tau^{\gamma}<T)}{P_{\theta}(\tau<T)}.\\\notag
\end{align}

\begin{align}\label{9}
&\int_{-E}^{\infty} \int_0^{\infty} \fmu e^{-\gamma t} e^{-\delta x} dt dx =
\frac{\vb}{\vb + \gamma}\frac{\int_{-E}^{\infty}   P_{\theta}( B^{\mu}(\tau^{\gamma})\in dx, \tau^{\gamma}<T) e^{-\delta x} dx }{P_{\theta}(\tau<T)} \\\notag
 &=\frac{\vb}{(\vb + \gamma)P_{\theta}(\tau<T)} \int_{-E}^{\infty} P_{\theta}\left(\underset{0<s<\tau^{\gamma}}{\inf} \Bsmu >-E , B^{\mu}(\tau^{\gamma})\in dx\right) e^{-\delta x} dx \\\notag
 &=
\frac{\vb}{(\vb + \gamma)P_{\theta}(\tau<T)} \int_{-E}^{\infty} \left[P_{\theta}( B^{\mu}(\tau^{\gamma})\in dx) e^{-\delta x}
- P_{\theta}\left(\underset{0<s<\tau^{\gamma}}{\inf} \Bsmu \leq -E , B^{\mu}(\tau^{\gamma})\in dx\right)e^{-\delta x} \right]dx  \\\notag
 &=
\frac{\vb + \gamma}{\sqrt{2(\vb+\gamma)+\Fmu^2}}
\frac{\vb}{(\vb + \gamma)P_{\theta}(\tau<T)}\cdot\\\notag
&\quad\cdot\int_{-E}^{\infty}
\left[e^{(\Fmu -\delta)x - |x| \sqrt{2(\vb+\gamma)+\Fmu^2} }-
e^{-x(\sqrt{2(\vb+\gamma)+\Fmu^2} +\delta-\Fmu)}e^{-2E\sqrt{2(\vb+\gamma)+\Fmu^2}}\right]dx  \\\notag
&= \frac{\vb}{\sqrt{2(\vb+\gamma)+\Fmu^2 }P_{\theta}(\tau<T)}\cdot \left[\int_{0}^{\infty}e^{-x(\sqrt{2(\vb+\gamma)+\Fmu^2} +\delta-\Fmu)}
 +\int_{-E}^{0}e^{x(\sqrt{2(\vb+\gamma)+\Fmu^2} -\delta+\Fmu)}\right.\\\notag
 &\qquad\qquad\qquad\qquad\qquad\qquad\qquad \left.- e^{ -2E\sqrt{2(\vb+\gamma)+\Fmu^2}} \int_{-E}^{\infty}e^{-x(\sqrt{2(\vb+\gamma)+\Fmu^2} +\delta-\Fmu)}    \right] \\\notag
 &=\frac{\vb(1-e^{-E(\sqrt{2(\vb+\gamma)+\Fmu^2} + \Fmu -\delta)})}{(\vb +\gamma+\delta\Fmu - \delta^2 / 2 )P_{\theta}(\tau<T)}.\\\notag
\end{align}
Where the first equality follows by ~(\ref{8}).

\textbf{Constructing the integral equation.} The DM chooses the risky arm, until the minimum between the stopping time of the first jump $\tau$ and the stopping time $T$. We distinguish between two cases.  In case the DM stops before the time of the first jump $\tau$, we calculate the expected discounted payoff from the risky arm until the stopping time $T$, and the expected discounted payoff from the safe arm after the stopping time $T$. In case the first jump occurs before the stopping time $T$, we calculate the expected discounted payoff received from the continuous part $\pita$ until time $\tau$. We add the expected discounted payoff from the first jump, and the expected discounted continuation payoff, updating both the posterior $p_{t,h,x}$, and the intercept $E+G_h +x$, according to the time of the first jump, the first jump's size, and the value of the continuous part of the payoff.\\\\
The notation used are as follows: $P_\theta (\tau>T)$ is the probability that the DM switches to the safe arm before a jump occurs. If $\tau>T$ then the expected payoff from the risky arm is $E_{\theta}[\int_0^{T}re^{-rt} d\pita | \tau > T]$, and the expected payoff from the safe arm is $E_{\theta}[ \int_T^{\infty} re^{-rt} sdt | \tau > T]$.
$P_\theta (\tau < T)$ is the probability a jump occurs before the DM switches to the safe arm. If $\tau<T$ then the expected payoff until the first jump is $E_{\theta}[\int_0^{\tau}re^{-rt} d\pita | \tau < T]$. $\fmu$ is the probability that the first jump occurs in the interval $[t,t+dt)$, and $\Bmutau$ belongs to $[x,x+dx)$, given that a jump occurs before the DM switches to the safe arm. $re^{-rt} \frac{1}{\bar{\nu}}\int h\nu (dh) $ is the expected discounted payoff from the first jump, and $\frac{1}{\bar{\nu}}\int\nu e^{-rt} U(p_{t,h,x},E+G_h+x))$ is the expected discounted continuation payoff, updating both the posterior $p_{t,h,x}$, and the intercept $E+G_h+x $ at time $t$. With this notation, the expected payoff is as follows: 
\begin{align}\label{integral equation}
.
\end{align}
$U(p_0,E) = p_0 \Patau \left[ E_{\ta}[\int_0^{\tau}re^{-rt} d\pita | \tau < T]
\right.$\\\\
$\left.+\int_{-E}^{\infty}\int_0^{\infty}  \fa  \left( re^{-rt} \frac{1}{\vg}\int \vh h  + \frac{1}{\vg}\int\vh e^{-rt} U(p_{t,h,x},E+G_h+x)\right)dtdx  \right] $\\\\
$
+p_0 \Patauu \left[ E_{\ta}[\int_0^{T}re^{-rt} d\pita | \tau > T] + E_{\ta}[ \int_T^{\infty} re^{-rt} sdt | \tau > T]\right] $\\\\
$
+(1-p_0) \Pbtau \left[ E_{\tb}[\int_0^{\tau}re^{-rt} d\pita | \tau < T]
 \right. $\\\\
$\left.+ \int_{-E}^{\infty}\int_0^{\infty}  \fb  \left( re^{-rt} \frac{1}{\vvg}\int \vvh h  + \frac{1}{\vvg}\int\vvh e^{-rt} U(p_{t,h,x},E+G_h+x)\right)dtdx  \right] $\\\\
$+(1-p_0) \Pbtauu \left[ E_{\tb}[\int_0^{T}re^{-rt} d\pita | \tau > T] + E_{\tb}[ \int_T^{\infty} re^{-rt} sdt | \tau > T]\right].
$\\\\\\

By Eqs.\ ~(\ref{3}), ~(\ref{4}) and ~(\ref{5}), this expression is equal to\\\\
$\ai p_0 (1-\Patauu) (1-E_{\ta}[e^{-r\tau}|\tau<T])$\\\\$
+  p_0 \Patau \Hi \int_{-E}^{\infty}\int_0^{\infty}  \fa  re^{-rt}dtdx$\\\\$
+  p_0 \Patau \int_{-E}^{\infty}\int_0^{\infty}  \fa   \frac{1}{\vg}\int e^{-rt} U(p_{t,h,x},E+G_h+x) \vh dtdx$\\\\$
+\ai p_0\Patauu (1-E_{\ta}[e^{-rT}| \tau > T]) +s  p_0\Patauu E_{\ta}[ e^{-rT}| \tau > T  ]$\\\\$
+\bi (1- p_0) \Pbtauu) (1-E_{\tb}[e^{-r\tau}|\tau<T])  $\\\\$
+ (1- p_0) \Pbtau \HHi \int_{-E}^{\infty}\int_0^{\infty}  \fb  re^{-rt}dtdx$\\\\$
+  (1- p_0)\Pbtau \int_{-E}^{\infty}\int_0^{\infty}  \fb   \frac{1}{\vvg}\int e^{-rt} U(p_{t,h,x},E+G_h+x)\vvh dtdx $\\\\$
+\bi (1- p_0)\Pbtauu (1-E_{\tb}[e^{-rT}| \tau > T] + s  (1- p_0)\Pbtauu E_{\tb}[ e^{-rT}| \tau > T  ].$\\\\\\
By Eqs.\ ~(\ref{6}), ~(\ref{7}) and ~(\ref{8}), this expression is equal to\\\\
$
\ai p_0 \Patauu\left(1-\frac{\Patauur}{\Patauu}\right)  +sp_0 \Patauu\frac{\Patauur}{\Patauu}$\\\\$
+p_0(1-\Patauu) \Hi r \cdot\frac{\vg}{\vg+r}\cdot\frac{1-\Patauur}{1-\Patauu}$\\\\$
+\ai p_0 (1-\Patauu) \left(1- \frac{\vg}{\vg+r}\cdot\frac{1-\Patauur}{1-\Patauu}\right)$\\\\$
+\bi (1- p_0) \Pbtauu\left(1-\frac{\Pbtauur}{\Pbtauu}\right) +s(1- p_0) \Pbtauu\frac{\Pbtauur}{\Pbtauu}$\\\\$
+(1- p_0)(1-\Pbtauu) \HHi r \cdot\frac{\vvg}{\vvg+r}\cdot\frac{1-\Pbtauur}{1-\Pbtauu}$\\\\$
+\bi (1- p_0) (1-\Pbtauu) \left(1- \frac{\vg}{\vg+r}\cdot\frac{1-\Pbtauur}{1-\Pbtauu}\right)$\\\\$
+
\int_{-E}^{\infty}\int_0^{\infty} \int  e^{-rt} U(p_{t,h,x},E+G_h+x) \cdot$\\\\$
\quad\cdot\left[p_0\Patau\fa\vh/\vg + (1-p_0)\Pbtau\fb\vvh/\vvg  \right]dtdx $\\\\\\
$=p_0 \left(\frac{\ai r}{\vg+r}+\frac{\vg \Hi r}{\vg+r} \right) + (1-p_0)\left(  \frac{\bi r}{\vvg+r} + \frac{\vvg \HHi r}{\vvg+r} \right) $\\\\$
+ p_0 \Patauur \left(s - \frac{\ai r}{\vg+r}  - \frac{\vg \Hi r}{\vg+r}  \right) + (1-p_0)\Pbtauur \left(s - \frac{\bi r}{\vvg+r}  - \frac{\vvg \HHi r}{\vvg+r}  \right)
$\\\\$
+\int_{-E}^{\infty}\int_0^{\infty} \int  e^{-rt} U(p_{t,h,x},E+G_h+x) \cdot$\\\\$
\cdot\left[p_0\Patau\fa\vh/\vg + (1-p_0)\Pbtau\fb\vvh/\vvg  \right] dtdx$\\\\\\
$=p_0 \left(\frac{\ai r}{\vg+r}+\frac{\vg \Hi r}{\vg+r} \right) + (1-p_0)\left(  \frac{\bi r}{\vvg+r} + \frac{\vvg \HHi r}{\vvg+r} \right) $\\\\$
+ p_0 e^{-E (F_{\ai} + \sqrt{2(\bar{\vi}+r)+F_{\ai}^2})} \left(s - \frac{\ai r}{\vg+r}  - \frac{\vg \Hi r}{\vg+r}  \right) $\\\\$
+ (1-p_0)e^{-E (F_{\bi} + \sqrt{2(\bar{\vvi}+r)+F_{\bi}^2})} \left(s - \frac{\bi r}{\vvg+r}  - \frac{\vvg \HHi r}{\vvg+r}  \right)
$\\\\$
+\int_{-E}^{\infty}\int_0^{\infty} \int  e^{-rt} U(p_{t,h,x},E+G_h+x) \cdot$\\\\$
\cdot\left[p_0\Patau\fa\vh/\vg + (1-p_0)\Pbtau\fb\vvh/\vvg  \right]dtdx .$\\\\
Simplifying the last expression, the integral equation is
\begin{align}\notag
U(p,E) &= Ap+B(1-p) +Cpe^{-m_1 E}+ D(1-p)e^{-m_2 E}\\\tag{IE}
&+ \int_{-E}^{\infty}\int_0^{\infty} \int  e^{-rt} U(p_{t,h,x},E+G_h+x) g(x,t,h)dtdx,\notag
\end{align}
where $A, B, C,$ and $D$ are constants, and\\
\begin{align}
g(x,t,h) &= \left[p_0\Patau\fa\vh/\vg\right. \\\notag
            &+  \left. (1-p_0)\Pbtau\fb\vvh/\vvg  \right] .
\end{align}

We show now that (IE) admits a unique solution $U$ in the region $[0,1]\times[0,\infty)$.\\\\
\textbf{Boundary values}\\
First, we find the values of $U(p,E)$ on the boundary of the region $[0,1]\times [0,\infty]$.
Note that $ U(p,0) = s,$ and $ U(p,\infty) = p\gi +(1-p)\ggi$ for every $p\in[0,1]$, .\\
We now argue that there is a unique solution for (IE) when $p=0$.
$U(0,E)$ is a function of $E$ with two boundary conditions, at $E=0$ and at $E=\infty$.
Suppose that $U(0,E), V(0,E)$ solve (IE). Then $W(0,E):=U(0,E)- V(0,E) $ satisfies $W(0,E) =  \int_{-E}^{\infty}\int_0^{\infty} \int  e^{-rt} W(0,E+G_h+x) g(x,t,h)dtdx$, and $W(0,0)=W(0,\infty)=0$. Let $\hat{E}$ be a critical point, where $W$ achieves its maximum.
Assume to the contrary that $W(0,\hat{E})>0$. By ~(\ref{8}) it follows that $\int_{-\hat{E}}^{\infty}\int_0^{\infty} \int  e^{-rt}  g(x,t,h) <1$. Therefore,
\begin{align}\notag
W(0,\hat{E}) &=  \int_{-\hat{E}}^{\infty}\int_0^{\infty} \int  e^{-rt} W(0,\hat{E}+G_h+x) g(x,t,h)\\\notag
             &\leq W(0,\hat{E})\int_{-\hat{E}}^{\infty}\int_0^{\infty} \int  e^{-rt}g(x,t,h)<W(0,\hat{E}),\notag
             \end{align}
%
which implies that $W(0,E)\leq 0$. Similarly, one can obtain that $W(0,E)\geq 0 $, so that $W(0,E)\equiv 0$, and the solution is unique.\\
Similar arguments show that (IE) admits a unique solution on $[0,\infty)$ when $p=1$.\\

Since $U(p,E)$ is uniquely determined on the boundary of the region $ [0,1]\times[0,\infty]$, similar arguments
show the uniqueness of the solution $[0,1]\times [0,\infty)$.
Using Eqs.\ ~(\ref{1}) - ~(\ref{9}) one can verify that the solution for (IE) is\\
$U(p_0,E) = p_0\gi + (1-p_0)\ggi + (s-\gi)p_0 e^{-(\at-\bt)(\alpha+1) E}+ (s-\ggi)(1-p_0)e^{-(\at-\bt)\alpha E}. $

By substituting $E : = \frac{\sigma}{\ai -\bi}\times \left[\ln\left(\frac{q_0}{1-q_0}\right)-\ln\left(\frac{p'}{1-p'}\right)\right]$, we get Eq.\ ~(\ref{E:UpE1}), as desired.
\end{proof}

\subsection{Information Pricing}

\begin{lem}\label{f1}
Let $f_{a,b}(\eta)$ be the function defined in (\ref{E:f1}). The equation $f_{a,b}(\eta)=0$ admits a unique solution in the interval $(0,\infty)$.
\end{lem}
The proof is similar to the proof of Lemma \ref{alpha}, and therefore omitted.
We turn to the proof of Theorem \ref{pricing1}, which is analogous to the proof of Theorems \ref{Up1} and \ref{UpE1}.
\begin{proof}[\textbf{Proof of Theorem \ref{pricing1}}]

Suppose that until time $t$, the DM chose the risky arm, and observed the jumps $h^a_1,...h^a_n$ (resp. $h^b_1,...h^b_m$) from the process $(X^a (t))$ (resp. $(X^b (t))$). Let $Y_B^a (t)$ (resp. $Y_B^b (t)$) be the Brownian motion with drift component of $(X^a (t))$ (resp. $(X^b (t))$) at time $t$. Note that $Y_B^j (t) \sim N(\mu^j t  , (\sigma^j)^2 t )$, $j\in\{ a, b \}$. Let $q_t: = P_t(\ta |h_1,...h_n;  \,  h^b_1,...h^b_m ;  \,  Y_B^a (t) ;  \,  Y_B^b (t);\,  q_0)$ be the posterior belief of the DM. The odd ratio of the posterior belief is 
\begin{align}\label{posterior2}
&\frac{q_t}{1-q_t} = \frac{q_0
\frac{1}{\sqrt{2\pi t}\sigma^a}
e^{-\frac{(Y_B^a (t)-\aia t)^2}{2(\sigma^a)^2 t}}
\frac{1}{\sqrt{2\pi t}\sigma^b}
e^{-\frac{(Y_B^b (t)-\aib t)^2}{2(\sigma^b)^2 t}}
e^{-\vga t}\prod_{t-}\vhia
e^{-\vgb t}\prod_{t-}\vhib }
{(1- q_0)
\frac{1}{\sqrt{2\pi t}(\sigma^a)}
e^{-\frac{(Y_B^a (t)-\bia t)^2}{2(\sigma^b)^2 t}}
\frac{1}{\sqrt{2\pi t}(\sigma^b)}
e^{-\frac{(Y_B^b (t)-\bia t)^2}{2(\sigma^b)^2 t}}
e^{-\vvga t}\prod_{t-}\vvhia
e^{-\vvgb t}\prod_{t-}\vvhib} \\\notag
&= \frac{q_0
e^{\aia Y_B^a (t)/(\sigma^a)^2 - (\aia)^2 t /2(\sigma^a)^2 }
e^{\aib Y_B^b (t)/(\sigma^b)^2 - (\aib)^2 t /2(\sigma^b)^2 }
e^{-\vga t} \prod_{t-} \vhia
e^{-\vgb t} \prod_{t-} \vhib}
{(1- q_0)
e^{\bia Y_B^a (t) /\sigma^2 - (\bia)^2 t /2(\sigma^a)^2 }
e^{\bib Y_B^b (t) /\sigma^2 - (\bib)^2 t /2(\sigma^b)^2 }
e^{-\vvga t} \prod_{t-} \vvhia
e^{-\vvgb t} \prod_{t-} \vvhib}.\notag
\end{align}

Indeed, $\frac{1}{\sqrt{2 \pi t}\sigma^a} e^{-\frac{(Y_B^a (t)-\aia t)^2}{2(\sigma^a)^2 t}}$
(resp. $\frac{1}{\sqrt{2 \pi t}\sigma^b} e^{-\frac{(Y_B^b (t)-\aib t)^2}{2(\sigma^b)^2 t}}$)
is the probability of observing $Y_B^a (t)$ (resp. $Y_B^b (t)$), given the type $\theta_i$,
and $e^{-\bar{\nu}^a_i t}\frac{(\bar{\nu}^a_i t)^n}{n!}\prod_{t-}\frac{\nu^a_i (dh^a_j)}{\bar{\nu}^a_i}  $
(resp. $e^{-\bar{\nu}^b_i t}\frac{(\bar{\nu}^b_i t)^m}{m!}\prod_{t-}\frac{\nu^b_i (dh^b_j)}{\bar{\nu}^b_i}  $ )
is the probability of receiving the $n$ (resp. $m$) jumps that occurred until time $t$ from $(X^a_i (t))$ (resp. $(X^b_i (t))$), given the type $\theta_i$. The first equality in ~(\ref{posterior2}) is the Bayesian belief updating, using the independence of the L\'{e}vy processes $X^a_i (t)$ and $X^b_i (t)$, and the independence of the components in the L\'{e}vy-Ito decomposition, given the type of the risky arm, and the second equality is obtained by eliminating common components. \\\\

Suppose the DM 
follows a cut-off strategy  $\kappa'$ with cut-off point $p'$.
If $q_0 \leq p'$, the DM will always choose the safe arm. If $q_0>p'$ the DM will initially choose the risky arm. 
Then the DM chooses the risky arm as long as $q_t > p'$, which, by Eq.\ ~(\ref{posterior}), is equivalent to
\begin{align}
&\dfrac{q_0}{1-q_0}
\times\dfrac{e^{\aia Y_B^a (t) /(\sigma^a)^2 - (\aia)^2 t /2(\sigma^a)^2 -\vga t}}
      {e^{\bia Y_B^a (t) /(\sigma^a)^2 - (\bia)^2 t /2(\sigma^a)^2 -\vvga t}}
\times\dfrac{e^{\aib Y_B^b (t) /(\sigma^b)^2 - (\aib)^2 t /2(\sigma^b)^2 -\vgb t}   }
      {e^{\bib Y_B^b (t) /(\sigma^b)^2 - (\bib)^2 t /2(\sigma^b)^2 -\vvgb t} }\\\notag
&\times\prod_{t-} \vhvva
\times\prod_{t-} \vhvvb
> \dfrac{p'}{1-p'}.
\end{align}
By taking the natural logarithm, and rearranging the resulting terms, we obtain that this inequality is equivalent to
%
\begin{align}\label{posteriorq2}
&\left(\frac{\aia-\bia}{\sigma^a}\right)Y_B^a (t)+
\left(\frac{\aib-\bib}{\sigma^b}\right)Y_B^b (t)>\\\notag
&\left(\frac{(\aia)^2-(\bia)^2}{2(\sigma^a)^2}+
\frac{(\aib)^2-(\bib)^2}{2(\sigma^b)^2}+ \vga+\vgb-\vvga-\vvgb\right)t\\\notag
&-\left( \ln\left(\frac{q_0}{1-q_0}\right)-\ln\left(\frac{p'}{1-p'}\right) \right)-
\sum_{t-} \lnvhvva-
\sum_{t-} \lnvhvvb.
\end{align}
%
Since $Y_B^j (t)=\mu^a + \sigma^a Z^a (t)$, $k\in\{a,b\}$, it follows that Eq.\ ~(\ref{posteriorq2}) is equivalent to
\begin{align}\label{Bmu}
B^{\muga,\mugb}(t):&=\frac{(\ata-\bta)Z^a (t)+(\atb-\btb)Z^b (t)}{\sqrt{(\ata-\bta)^2 +(\atb-\btb)^2}} +\left[\frac{-(\vga+\vgb-\vvga-\vvgb)
t}{\sqrt{(\ata-\bta)^2 +(\atb-\btb)^2}}\right.\\\notag
&\left.+
\frac{2\tilde{\mu}^a (\ata-\bta)-((\ata)^2-(\bta)^2)
+2\tilde{\mu}^b (\atb-\btb)-((\atb)^2-(\btb)^2)}{2\sqrt{(\ata-\bta)^2 +(\atb-\btb)^2}}\right]t\\\notag
&>\frac{-\left(\ln\left(\frac{q_0}{1-q_0}\right)-\ln\left(\frac{p'}{1-p'}\right)
\right)
-\sum_{t-} \lnvhvva-
\sum_{t-} \lnvhvvb
}{\sqrt{(\ata-\bta)^2 +(\atb-\btb)^2}},
\end{align}
where $\tilde{\mu}^j_i=\frac{\mu^j_i}{\sigma^j}$, for $j\in\{a,b,c\}$, and $i\in\{ \high ,\low \}$.
Notice that $B^{\muga,\mugb}(t)$ is a standard Brownian motion with drift:
\begin{align}\notag
F^{\muga,\mugb}:&=
\left[\frac{-(\vga+\vgb-\vvga-\vvgb)
t}{\sqrt{(\ata-\bta)^2 +(\atb-\btb)^2}}\right.\\\notag
&\left.+\frac{2\tilde{\mu}^a (\ata-\bta)-((\ata)^2-(\bta)^2)
+2\tilde{\mu}^b (\atb-\btb)-((\atb)^2-(\btb)^2)}{2\sqrt{(\ata-\bta)^2 +(\atb-\btb)^2}}\right].
\end{align}

We construct an integral equation similar to the one in the proof of Theorem \ref{UpE1}.\\

Let $\tau$ be the stopping time of the first jump from the process $(X^a (t)+ X^b (t))$. Let the stopping time $T$ be the first time $t$ that satisfies
\begin{align}
B^{\muga,\mugb}(t)\leq
\frac{-\left(\ln\left(\frac{q_0}{1-q_0}\right)-\ln\left(\frac{p'}{1-p'}\right)
\right)
-\sum_{t-} \lnvhvva-
\sum_{t-} \lnvhvvb
}{\sqrt{(\ata-\bta)^2 +(\atb-\btb)^2}}.
\end{align}

Denote $G^a_h:=\frac{\lnvhvva}{\sqrt{(\ata-\bta)^2 +(\atb-\btb)^2}}$ (resp. $G^b_h:=\frac{ \lnvhvvb}{\sqrt{(\ata-\bta)^2 +(\atb-\btb)^2}}$) the contribution of a jump of size $h^a$ (resp. $h^a$) received from the process $X^a (t)$ (resp. $X^a (t)$), and
$E:= \frac{-\left(\ln\left(\frac{q_0}{1-q_0}\right)-\ln\left(\frac{p'}{1-p'}\right)
\right)
}{\sqrt{(\ata-\bta)^2 +(\atb-\btb)^2}}$, the intercept of the right-hand side of Eq.\ ~(\ref{Bmu}) at $t=0$.


\textbf{Constructing the integral equation.} The DM chooses the risky arm, until the minimum between the stopping time of the first jump $\tau$ and the stopping time $T$. We distinguish between two cases. In case the DM stops before the time of the first jump $\tau$, we calculate the expected discounted payoff from the process $(X^a (t)+ X^c (t))$ until the stopping time $T$, and the expected discounted payoff from the safe arm after the stopping time $T$. In case the first jump occurs before the stopping time $T$, we calculate the expected discounted payoff received from the process $(X^a (t) + X^c (t))$ until time $\tau$.
If the first jump was received from the process $(X^a (t))$,
we add the expected discounted payoff from the first jump, and the expected discounted continuation payoff, updating both the posterior $p^a_{t,h^a,x}$, and the intercept $E+G^a_h +x$, according to the time of the first jump, the first jump's size, and the value of the continuous part of the process $(X^a (t)+ X^b (t))$; while if the first jump was received from the process $(X^b (t))$,
we add only the expected discounted continuation payoff, updating both the posterior $p^b_{t,h^b,x}$, and the intercept $E+G^b_h +x$, according to the time of the first jump, the first jump's size, and the value of the continuous part of the process $(X^a (t)+ X^b (t))$.


The posterior $p^j_{t,h^j,x}$ is updated as follows:
\begin{align}\notag
&p^j_{t,h^j,x} := P(\ta | \tau<T, \tau =t, B^{\ata,\bta}(\tau) \in dx, h^j ) \\\notag
       &= \frac{ P( \tau<T, \tau =t, B^{\muga,\mugb}(\tau) \in dx, h^j|\ta )P(\ta)}
               {P( \tau<T, \tau =t, B^{\muga,\mugb}(\tau) \in dx, h^j|\ta )P(\ta) +
                P( \tau<T, \tau =t, B^{\muga,\mugb}(\tau) \in dx, h^j|\tb )P(\tb)}  \\\notag
       &=\frac{q_0 \Patau \faj \frac{\nu^j_{\high}(dh^j)}{\vgj}}{q_0 \Patau \faj \frac{\nu^j_{\high}(dh^j)}{\vgj} + (1-q_0) \Pbtau \fbj \frac{\nu^j_{\low} (dh^j)}{\vvgj}}.    \notag
\end{align}

%
%
%
%
The integral equation is as follows:
\begin{align}\label{IE3}
&U(p_0,E) = p_0 \Patau \left[ E_{\ta}\left[\left.\int_0^{\tau}re^{-rt} dY^a_B(t) \right| \tau < T\right]
\right.+E_{\ta}\left[\left.\int_0^{\tau}re^{-rt} dY^c (t) \right| \tau < T\right] \\\notag
&+\int_{-E}^{\infty}\int_0^{\infty}  \faj   re^{-rt} \frac{\vga\Hia}{\vga+\vgb}dtdx
\\\notag
&+\int_{-E}^{\infty}\int_0^{\infty}  \faj  e^{-rt}\int_h U(p^a_{t,h^a,x},E+G^a_h+x)\left(\frac{\vha}{\vga+\vgb}\right)dtdx
\\\notag
&\left.+\int_{-E}^{\infty}\int_0^{\infty}  \faj  e^{-rt}\int_h U(p^b_{t,h^b,x},E+G^b_h+x)\left(\frac{\vhb}{\vga+\vgb}\right)dtdx\right]
\\\notag
&+p_0 \Patauu \left[ E_{\ta}\left[\left.\int_0^{T}re^{-rt} dY^a_B(t) \right| \tau > T\right] +E_{\ta}\left[\left.\int_0^{T}re^{-rt} dY^c(t) \right| \tau > T\right]\right.\\\notag
&\qquad\qquad\qquad\quad\left.+ E_{\ta}\left[\left. \int_T^{\infty} re^{-rt} sdt \right| \tau > T\right]\right] \\\notag
&
+(1-p_0) \Pbtau \left[ E_{\tb}\left[\left.\int_0^{\tau}re^{-rt} dY^a_B(t) \right| \tau < T\right]
 \right.+E_{\tb}\left[\left.\int_0^{\tau}re^{-rt} dY^c (t) \right| \tau < T\right]  \\\notag
&+\int_{-E}^{\infty}\int_0^{\infty}  \fbj   re^{-rt} \frac{\vvga\HHia}{\vvga+\vvgb}dtdx
\\\notag
&+\int_{-E}^{\infty}\int_0^{\infty}  \fbj  e^{-rt}\int_h U(p^a_{t,h^a,x},E+G^a_h+x)\left(\frac{\vvha}{\vvga+\vvgb}\right)dtdx
\\\notag
&\left.+\int_{-E}^{\infty}\int_0^{\infty}  \fbj  e^{-rt}\int_h U(p^b_{t,h^b,x},E+G^b_h+x)\left(\frac{\vvhb}{\vvga+\vvgb}\right)dtdx\right]
\\\notag
&+(1-p_0) \Pbtauu \left[  E_{\tb}\left[\left.\int_0^{T}re^{-rt} dY^a_B(t) \right| \tau > T\right] +E_{\tb}\left[\left.\int_0^{T}re^{-rt} dY^c(t) \right| \tau > T\right]\right.\\\notag
&\qquad\qquad\qquad\qquad\quad\;\left.+ E_{\tb}\left[\left. \int_T^{\infty} re^{-rt} sdt \right| \tau > T\right]\right].
\end{align}
Similarly to Eqs.\ ~(\ref{4}) and ~(\ref{5}), using the stochastic integral properties for L\'{e}vy processes, it follows that for $i\in\{\high,\low\}$
\begin{equation}\notag
E_{\theta_i}\left[\left.\int_0^{T}re^{-rt} dY^c(t) \right| \tau > T\right] = g^c_i \left(1-E_{\theta}[e^{-rT}|\tau>T] \right),
\end{equation} and
\begin{equation}\notag
E_{\theta_i}\left[\left.\int_0^{\tau}re^{-rt} dY^c(t) \right| \tau < T\right] = g^c_i \left(1-E_{\theta}[e^{-r\tau}|\tau<T] \right).
\end{equation}
By these equations, and Eqs.\ ~(\ref{2})- ~(\ref{9}), one can verify that for every $q_0>p'$, the unique solution of Eq.\ ~(\ref{IE3}) is

\begin{align}\notag
U(p_0,E)&=(s-\gia-\gic)p_0e^{-(\beta+1)\sqrt{(\ata-\bta)^2 +(\atb-\btb)^2}E} \\\notag
& +(s-\ggia-\ggic)(1-p_0)e^{-\beta\sqrt{(\ata-\bta)^2 +(\atb-\btb)^2}E}\\\notag
&+p_0(\gia+\gic)+(1-p_0)(\ggia+\ggic). \notag
\end{align}
Therefore, for every $q_0>p'$, the expected discounted payoff is
\begin{align}\notag
V_{p'} (p_0,q_0)&=  (s-\gia-\gic)p_0\left(\frac{1-q_0}{q_0}\right)^{\beta+1} \left(\frac{p'}{1-p'}\right)^{\beta +1} \\\notag
& +(s-\ggia-\ggic)(1-p_0)\left(\frac{1-q_0}{q_0}\right)^{\beta} \left(\frac{p'}{1-p'}\right)^{\beta }\\\notag
&+p_0(\gia+\gic)+(1-p_0)(\ggia+\ggic), \notag
\end{align}
as desired.

Now, if $q_0=p_0$, using the same method to that used to prove Theorem \ref{Up1}, one can show that the optimal cut-off strategy is
\[
 p^*= \frac{\beta (s- \ggia - \ggic)}{(\beta + 1)(\gia + \gic-s) + \beta(s - \ggia - \ggic)}.
\]

\end{proof}

\begin{proof}[\textbf{Proof of Corollary \ref{corollary}}]

Let $\beta_a$ be the unique solution of
\begin{align}\label{fa}
f_{a} (\eta):&=
\int\vvha\left(\frac{\vvha}{\vha}\right)^{\eta}
+ \eta (\vga-\vvga)
- \vvga
+\frac{1}{2}(\eta+1)\eta \left(\frac{\aia-\bia}{\sigma^a}\right)^2 - r = 0\notag
\end{align} in $(0,\infty)$, and let $\beta_{a,b}$ be the unique solution of Eq.\ ~(\ref{E:f1}). Using similar arguments to those used to prove Lemma \ref{alpha}, one can show that $\beta_{a,b} < \beta_{a}$.
Let $p^*_{a,c} = \frac{\beta_{a} (s- \ggia - \ggic)}{(\beta_{a} + 1)(\gia + \gic-s) + \beta_{a}(s - \ggia - \ggic)}$ be the optimal cut-off of DM2. By Theorem \ref{pricing1}, the optimal expected payoff of DM2 is given by
\begin{align}\notag
U_{DM2}(p_0)&=V^{a,c}_{p^*_{a,c}} (p_0,p_0)\\\notag
&=\begin{cases}
s               &\text{if $p_0 \leq p^*_{a,c}$}, \\
p_0(\gia+\gic) +(1-p_0)(\ggia+\ggic) +C_{\beta_{a}} (1-p_0)(\frac{1-p_0}{p_0})^{\beta_a}   &\text{if $p_0 > p^*_{a,c}$},
\end{cases}
\end{align}
where $C_{\beta_{a}} = \frac{s-\ggia-\ggic - p^*_{a,c} (\gia-\gic -\ggia-\ggic)}{(1-p^*_{a,c}) \left( \frac{1-p^*_{a,c}}{p^*_{a,c}}\right)^{\beta_{a}}}.$

By Theorem \ref{pricing1}, the expected payoff of DM1, using the same cut-off point $p^*_{a,c}$ is given by
\begin{align}\notag
V^{a,b,c}_{p^*_{a,c}} (p_0,p_0)=
\begin{cases}
s               &\text{if $p_0 \leq p^*_{a,c}$}, \\
p_0(\gia+\gic) +(1-p_0)(\ggia+\ggic) +C_{\beta_{a,b}} (1-p_0)(\frac{1-p_0}{p_0})^{\beta_{a,b}}   &\text{if $p_0 > p^*_{a,c}$},
\end{cases}
\end{align}
where $C_{\beta_{a,b}} = \frac{s-\ggia-\ggic - p^*_{a,c} (\gia-\gic -\ggia-\ggic)}{(1-p^*_{a,c}) \left( \frac{1-p^*_{a,c}}{p^*_{a,c}}\right)^{\beta_{a,b}}}.$

Since $\beta_{a,b} < \beta_{a}$, it follows that if $p_0 > p^*_{a,c}$, then $U_{DM2}(p_0)=V^{a,c}_{p^*_{a,c}} (p_0,p_0)< V^{a,b,c}_{p^*_{a,c}}  (p_0,p_0)\leq U_{DM1}(p_0)$, as desired.

\end{proof}

\subsection{Optimism vs Pessimism}\label{ss:optimism}
\begin{proof}[\textbf{Proof of Theorem \ref{optpes}}]

Substituting $q_0=p_0+\epsilon$ in Eq.\ ~(\ref{E:UpE1}) yields


\begin{align}\notag
V_{p_0} (\epsilon)
&=p_0 \gi+(1-p_0)\ggi + p_0 (s -\gi) \left(\fracp1*\right)^{\alpha+1} \left(\fraceps\right)^{\alpha +1}\\\notag
&+(1-p_0) (s-\ggi) \left(\fracp1*\right)^{\alpha} \left(\fraceps\right)^{\alpha },\notag
\end{align}
where $\alpha$ is the unique solution of ~(\ref{eta})
in $(0,\infty)$.\\\\
Simple algebraic manipulations yield
\begin{eqnarray}\notag
V_{p_0} (\epsilon)
&=&p_0\gi + (1-p_0)\ggi+ p_0 (s -\gi) \left(\fracp1*\right)^{\alpha+1} \left(\fraceps\right)^{\alpha +1} \\\notag
&&+(1-p_0) (s-\ggi) \left(\fracp1*\right)^{\alpha} \left(\fraceps\right)^{\alpha}\\\notag
&=& [p_0\gi +  (1-p_0)\ggi] + \left(\fracp1*\right)^{\alpha}\cdot\\\notag
 &&\cdot \left[  p_0 (s-\gi)\fracp1*\left(\fraceps\right)^{\alpha +1} +(1-p_0) (s-\ggi)\left(\fraceps\right)^{\alpha }\right]\\\notag
&=&[p_0\gi + (1-p_0)\ggi]+  \left(\frac{\alpha}{\alpha+1}\right) ^{\alpha} \left(\frac{s-\ggi}{\gi-s}\right) ^{\alpha}
\cdot \\\notag
&&\cdot\left[ - p_0 (s-\ggi)\frac{\alpha}{\alpha+1}\left(\fraceps\right)^{\alpha +1} +(1-p_0) (s-\ggi)\left(\fraceps\right)^{\alpha }\right]\\\notag
&=& [p_0\gi + (1-p_0)\ggi]+ \left(\frac{\alpha}{\alpha+1}\right) ^{\alpha} \left(\frac{s-\ggi}{\gi-s}\right) ^{\alpha}\frac{1}{p_0 (s-\ggi)}
\cdot \\\notag
&&\cdot \left[ - \frac{\alpha}{\alpha+1}\left(\fraceps\right)^{\alpha +1} +\frac{1-p_0}{p_0} \left(\fraceps\right)^{\alpha }\right].\\\notag
\end{eqnarray}
%
%
%
%
For every $ x \geq 0$ define $W(x)=V_p (x)-V_p(-x) $.
This is the difference between the payoff of an optimist and the payoff of a pessimist. Straightforward calculations show that\\
\[W'(x) = \frac{\alpha x}{p}\left[
-\left( \frac{1-(p+x)}{p+x}    \right)^{\alpha}\frac{1}{(p+x)^2 (1-(p+x))}
+\left( \frac{1-(p-x)}{p-x}    \right)^{\alpha}\frac{1}{(p-x)^2 (1-(p-x))}
   \right],
\]\\ so that $W(0)=W'(0)=0$.
Suppose $\alpha>1$. 
Since
$\frac{1}{(p+x)^3} < \frac{1}{(p-x)^3} $, and
$\left( \frac{1-(p+x)}{p+x}    \right)^{\alpha-1} < \left(\frac {1-(p-x)}{p-x}\right)    ^{\alpha-1}  $,
we get $W'(x)>0$ for every $x >0$ such that $\p1* < p\pm x \leq 1$, and so $W(x)>0$; an optimist will fare better.\\\\
If $0 < \alpha < 1$, it is easy to verify that $\p1* < \frac{\alpha+2}{3} $.
Since the function $\left(\frac{1-y}{y}    \right)^{\alpha}\frac{1}{y^2 (1-y)}$ decreases for $0<y<\frac{\alpha+2}{3}$
and increases for $\frac{\alpha+2}{3}<y<1$, it follows that for every $x >0$ such that $\p1* < p\pm x \leq \frac{\alpha+2}{3}$ we get
$W'(x)>0$; an optimist will fare better, and for every $x >0$ such that $  \frac{\alpha+2}{3} <p\pm x$ we get $W'(x)<0$, and so $W(x)<0$; a pessimist will fare better.
\end{proof}

\section{Figures}
In Figure 1 we depict a generic path of the L\'{e}vy payoff process $(Y(t))$. The contribution of the compound Poisson process (the jumps) appears in Figure 2. The posterior belief given the observations appears in Figure 3. Finally, Figure 4 shows the time-dependent cut-off, as well as the continuous part of the process (the Brownian motion with drift).
\begin{figure}[b]
\centering
    \includegraphics[scale=0.48]{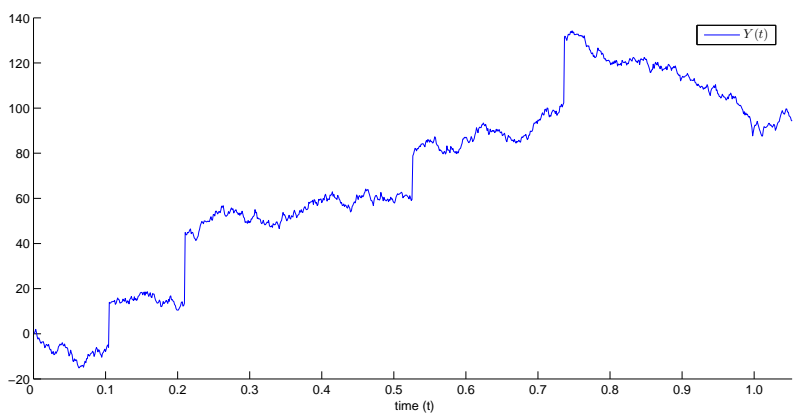}
    \caption{The payoff process }
\end{figure}

\begin{figure}[b]
\centering
 \begin{center}
    \includegraphics[scale=0.48]{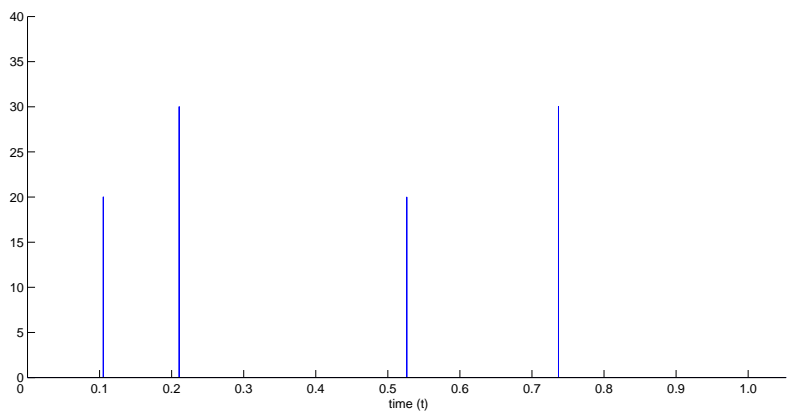}
    \caption{The Poisson arrivals}
\end{center}
\end{figure}

\begin{figure}[]
\centering
 \begin{center}
    \includegraphics[scale=0.52]{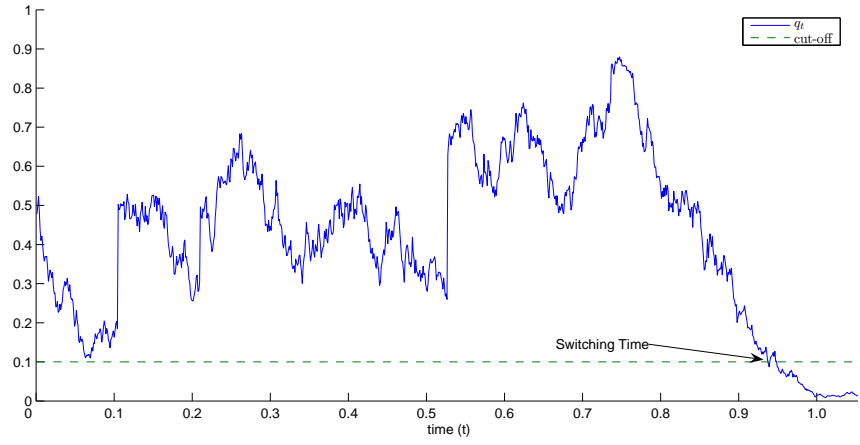}
    \caption{The posterior process}
\end{center}
\end{figure}

\begin{figure}[]
\centering
    \includegraphics[scale=0.52]{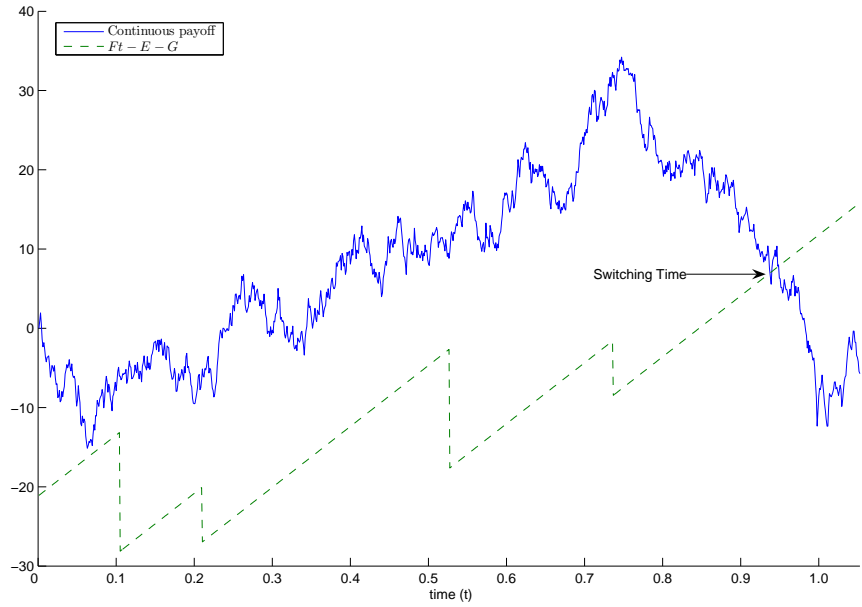}
    \caption{The strategy description }
\end{figure}
\newpage


\bibliographystyle{amsplain}
\bibliography{}

\end{document}